\documentclass[10pt,a4paper]{amsart}
\usepackage{verbatim}

\usepackage[toc]{appendix}
\usepackage[T1]{fontenc}

\usepackage{graphicx}
\usepackage{enumerate}
\usepackage{amsmath,amsfonts,amssymb}
\usepackage{color}
\def\loc{\operatorname{loc}}
\usepackage{cite}
\definecolor{citation}{rgb}{0.11,0.67,0.84}
\definecolor{formula}{rgb}{0.1,0.2,0.6}
\definecolor{url}{rgb}{0.11,0.67,0.84}
\usepackage{pgf,tikz}
\usepackage{mathrsfs}
\usepackage{fancyhdr}
\usepackage{dutchcal}

\usepackage{dsfont}

\newcommand{\reqnomode}{\tagsleft@false}

\usepackage[colorlinks,pdfpagelabels,pdfstartview = FitH,bookmarksopen = true,bookmarksnumbered = true,urlcolor=url,linkcolor = formula,plainpages = false,hypertexnames = false,citecolor = citation] {hyperref}

\def\dx{\,{\rm d}x}

\def\dy{\,{\rm d}y}

\def \d{\,{\rm d}}

\def\dist{\,{\rm dist}}

\allowdisplaybreaks
\makeatletter
\DeclareRobustCommand*{\bfseries}{%
  \not@math@alphabet\bfseries\mathbf
  \fontseries\bfdefault\selectfont
  \boldmath
}

\DeclareMathOperator*{\osc}{osc}

\makeatother

\newlength{\defbaselineskip}
\newcommand{\setlinespacing}[1]
           {\setlength{\baselineskip}{#1 \defbaselineskip}}
\newcommand{\mint}{\mathop{\int\hskip -1,05em -\, \!\!\!}\nolimits}

\newtheorem{theorem}{Theorem}

\newtheorem{definition}{Definition}
\newtheorem{remark}{Remark}[section]
\newtheorem{lemma}{Lemma}[section]

\numberwithin{equation}{section}
\allowdisplaybreaks[4]

\newcommand{\kk}{\kappa}

\def\en{\mathbb N}
\def\er{\mathbb R}

\newcommand{\ti}[1]{\tilde{#1}}

\newcommand\eps\varepsilon

\def\eqn#1$$#2$${\begin{equation}\label#1#2\end{equation}}

\newcommand{\be}{\begin{equation}}
\newcommand{\ee}{\end{equation}}

\newcommand{\rr}{\varrho}

\newcommand{\snr}[1]{\lvert #1\rvert}

\newcommand{\nr}[1]{\lVert #1 \rVert}

\newcommand{\N}{\mathbb{N}}

\def\name[#1, #2]{#1 #2}

\title[Multi-phase integrals]{Optimal gradient estimates for multi-phase integrals}

\author[De Filippis]{Cristiana De Filippis}  \address{Cristiana De Filippis\\Dipartimento di Matematica "Giuseppe Peano", Universit\`a di Torino\\ Via Carlo Alberto 10, 10123 Torino, Italy} \email{\url{cristiana.defilippis@unito.it}}

\begin{document}

\subjclass[2010]{35J20, 35J50, 35J60, 35J70, 49N60\vspace{1mm}} 

\keywords{Regularity, nonuniform ellipticity\vspace{1mm}}

\thanks{{\it Acknowledgements.}\ This work is supported by the University of Turin via the project "Regolarit\'a e propriet\'a qualitative delle soluzioni di equazioni alle derivate parziali".
\vspace{1mm}}

\maketitle

\begin{abstract}
We prove sharp reverse H\"older inequalities for minima of multi-phase variational integrals and apply them to Calder\'on-Zygmund estimates for nonhomogeneous problems.
\end{abstract}
\vspace{3mm}
\tableofcontents

\setlinespacing{1.08}
\section{Introduction}\label{si}
In this paper we complete the regularity theory started in \cite{deoh} for local minimizers of multi-phase functionals, i.e. variational integrals of the type
$$
W^{1,p}(\Omega)\ni w\mapsto \mathcal{H}(w,\Omega):=\int_{\Omega}\left[\snr{Dw}^{p}+\sum_{\nu=1}^{\kk}a_{\nu}(x)\snr{Dw}^{p_{\nu}}\right] \ \dx,
$$
where the modulating coefficients $\{a_{\nu}\}_{\nu\in 1}^{\kk}$ and exponents $(p,p_{1},\cdots,p_{\kk})$ satisfy
\begin{eqnarray}\label{ab}
0\le a_{\nu}(\cdot)\in C^{0,\alpha_{\nu}}(\Omega)\qquad \mbox{and}\qquad 1<p<\min_{\nu\in \texttt{I}_{\kk}}p_{\nu}
\end{eqnarray}
and it is $\texttt{I}_{\kk}:=\{1,\cdots,\kk\}$. Exponents $p$, $p_{\nu}$, $\alpha_{\nu}$ are related by the constraint
\begin{eqnarray}\label{bounds}
\frac{p_{\nu}}{p}\le 1+\frac{\alpha_{\nu}}{n}\qquad \mbox{for all} \ \ \nu\in\texttt{I}_{\kk},
\end{eqnarray}
which is sharp in the light of the counterexamples in \cite{bds,eslemi,fomami,Z1}. Precisely, our first achievement concerns some reverse H\"older type inequalities in the spirit of those obtained in \cite{comib,comi} for double phase problems.
\begin{theorem}\label{revh}
Under assumptions \eqref{ab}-\eqref{bounds}, let $v\in W^{1,p}_{\loc}(\Omega)$ be a local minimizer of functional $\mathcal{H}(\cdot)$ and $B_{\rr}(x_{0})\Subset \Omega$ be any ball with radius $\rr\in (0,1]$ and $J\ge 4$ be a constant. Then
\begin{itemize}
    \item in the degenerate regime $\texttt{deg}_{\textnormal{J}}(B_{\rr}(x_{0}))$ for all $d\ge 1$ it holds that
    \begin{eqnarray}\label{fh}
    \left(\mint_{B_{\rr/2}(x_{0})}\snr{Dv}^{d} \ \dx\right)^{1/d}&\le& cJ^{\Gamma}\left(\mint_{B_{\rr}(x_{0})}H(x,Dv) \ \dx\right)^{1/p},
    \end{eqnarray}
    with $c\equiv c(\texttt{data},\nr{H(\cdot,Dv)}_{L^{1+\delta_{g}}(B_{\rr}(x_{0}))},d)$ and $\Gamma\equiv \Gamma(\texttt{data}_{0})$;
    \item in the nondegenerate regime $\texttt{ndeg}_{\textnormal{J}}(B_{\rr}(x_{0}))$ or in the mixed one $\texttt{mix}_{\textnormal{J}}(B_{\rr}(x_{0}))$, for all $d\ge 1$, $\mu\in (0,1]$ it is
    \begin{eqnarray}\label{hf}
    \left(\mint_{B_{\rr/2}(x_{0})}\snr{Dv}^{d} \ \dx\right)^{1/d}&\le& c\rr^{-\mu}\left(\mint_{B_{\rr}(x_{0})}H(x,Dv) \ \dx\right)^{1/p},
    \end{eqnarray}
    for $c\equiv c(\texttt{data},\texttt{A},\nr{H(\cdot,Du)}_{L^{1+\delta_{g}}(B_{\rr}(x_{0}))},\mu,d)$.
\end{itemize}
\end{theorem}
We refer to Sections \ref{not} and \ref{crh} for more details on the terminology adopted in the above statement. A result analogous to the one described in Theorem \ref{revh} has been obtained in \cite[Theorem 4.1]{bbo} for generalized \cite{bb01,boh} triple phase problems, which in principle include also our functional $\mathcal{H}(\cdot)$. However, in \cite{bbo} to prove estimates similar to \eqref{fh}-\eqref{hf}, extra technical assumptions on $\{\alpha_{\nu}\}_{\nu=1}^{\kk}$ are required, i.e.: 
\eqn{ns}
$$
\max_{\nu\in \texttt{I}_{\kk}}\alpha_{\nu}\le 2\min_{\nu\in \texttt{I}_{\kk}}\alpha_{\nu},
$$
cf. \cite[(1.17), (1.22) and (6.8)]{bbo}. Condition \eqref{ns} seems to be unavoidable according to the arguments developed in \cite{bbo}, inspired by \cite{comi} and essentially relying on a boost of integrability that results from a combination of a Caccioppoli type inequality with the classical fractional Sobolev embedding theorem. In sharp contrast with what happens in \cite{comi}, the rate of nonhomogeneity in multi-phase problems is too high and causes competition among the H\"older continuity exponents $\{\alpha_{\nu}\}_{\nu=1}^{\kk}$. This drastically affects the integrability improvement granted by Sobolev embedding theorem and possibly leads to violations of the bounds in \eqref{bounds}. Here, we rather follow the approach of \cite{comib}, replace fractional Sobolev embedding theorem with a suitable fractional Gagliardo-Nirenberg inequality \cite{BrMi}, which matches the controlled gradient fractional differentiability assured by Caccioppoli inequality with the Morrey type result obtained in \cite[Theorem 2]{deoh}. Precisely, the idea consists in exploiting Gagliardo-Nirenberg inequality to translate the $\beta_{0}$-H\"older continuity of minima for arbitrary $\beta_{0}\in (0,1)$ consequence of \cite[Theorem 2]{deoh}, into gradient higher integrability up to any finite exponent, thus bypassing all structural obstructions due to the coexistence of multiple phases. In the light of \cite[Theorem 1]{deoh}, inequalities \eqref{fh}-\eqref{hf} do not add any substantial information on the regularity of minima of functional $\mathcal{H}(\cdot)$. Anyway, they turn out to be fundamental for instance when such minimizers play the role of comparison map in nonhomogeneous problems of the type
\begin{eqnarray}\label{czp}
W^{1,p}(\Omega)\ni w\mapsto \mathcal{G}(w,\Omega):= \int_{\Omega}\left[H(x,Dw)-\langle G(x,F),Dw\rangle\right] \ \dx
\end{eqnarray}
where $G\colon \Omega\times \mathbb{R}^{n}\to \mathbb{R}^{n}$ is a Carath\'eodory vector field so that
\begin{eqnarray}\label{assg}
\snr{G(x,z)}\le \Lambda \frac{H(x,z)}{\snr{z}}\qquad \mbox{for all} \ \ (x,z)\in \Omega\times \mathbb{R}^{n} \ \ \mbox{and some} \ \ \Lambda>0
\end{eqnarray}
and $F\colon \Omega \to \mathbb{R}^{n}$ verifies
\begin{eqnarray}\label{assf}
H(\cdot,F)\in W^{1,\gamma}_{\loc}(\Omega,\mathbb{R}^{n}) \qquad \mbox{with} \ \ \gamma>1.
\end{eqnarray}
For local minima of the functional in \eqref{czp} we have the following Calder\'on-Zygmund type result.
\begin{theorem}\label{czt}
Under assumptions \eqref{ab}, \eqref{bounds}, \eqref{assg}, \eqref{assf}, let $u\in W^{1,p}(\Omega)$ be a local minimizer of functional $\mathcal{G}(\cdot)$. Then the sharp Calrder\'on-Zygmund implication
\begin{eqnarray*}
H(\cdot,F)\in L^{\gamma}_{\loc}(\Omega) \ \Longrightarrow \ H(\cdot,Du)\in L^{\gamma}_{\loc}(\Omega) 
\end{eqnarray*}
holds for all $\gamma>1$. Moreover, for open sets $\Omega_{0}\Subset \tilde{\Omega}_{0}\Subset \Omega$ so that $\dist(\Omega_{0},\partial \ti{\Omega}_{0})\sim\dist(\ti{\Omega}_{0},\partial \Omega)\sim\dist(\Omega_{0},\partial \Omega)$. For every $\gamma>1$ there exists a radius $r_{*}>0$ and a constant $c\ge 1$, both depending on $(\texttt{data}_{\textnormal{cz}})$ such that
\begin{eqnarray}\label{czz}
\left(\mint_{B_{\rr/2}(x_{0})}H(x,Du)^{\gamma} \ \dx\right)^{1/\gamma}&\le&c\mint_{B_{\rr}(x_{0})}H(x,Du) \ \dx\nonumber \\
&+&c\left(\mint_{B_{\rr}(x_{0})}H(x,F)^{\gamma} \ \dx\right)^{1/\gamma},
\end{eqnarray}
for all balls $B_{\rr}(x_{0})\Subset \Omega_{0}$ with $\rr\in (0,r_{*})$.
\end{theorem}
We remark that Theorem \ref{czt} is not included in \cite[Theorem 1.1]{bbo} as we do not assume \eqref{ns}. Let us put our results into the context of the available literature. Multi-phase functionals provide the natural generalization of the double-phase energy
\begin{flalign*}
&W^{1,p}_{\loc}(\Omega)\ni w\mapsto \mathcal{P}(w,\Omega):=\int_{\Omega}\left[\snr{Dw}^{p}+a(x)\snr{Dw}^{q}\right] \ \dx,\nonumber \\
&\qquad \quad 0\le a(\cdot)\in C^{0,\alpha}(\Omega),\qquad \frac{q}{p}\le 1+\frac{\alpha}{n},
\end{flalign*}
first studied in \cite{Z1,Z2}, with emphasis about homogeneization and on the possible occurrence of Lavrentiev phenomenon and later on, regularity has been obtained in \cite{bcm2,bacomi,comib,comi}, see also \cite{bh1,comicz,demicz} concerning Calder\'on-Zygmund estimates, \cite{demima} on the general vectorial setting and the manifold-constrained case, \cite{chde} about potential theoretic considerations and \cite{ciccio} for sharp regularity of nonhomogenous systems with double phase structure and related obstacle problems and \cite{bb90,bacomi1,bb01,bbo,boh,by1,deoh,frzz,ok,rata} for further extensions and more general models. The peculiarity of the double phase energy is the subtle interaction between the $p$-phase and the $(p,q)$-phase, whose alternance is controlled by the modulating coefficient $a(\cdot)$: in proximity of the zero level set $\left\{x\in \Omega\colon a(x)=0\right\}$, the integrand in $\mathcal{P}(\cdot)$ behaves as the $p$-Laplacian, while in correspondence of the positivity set of $a(\cdot)$ it acts as a $\Delta_{2}$-Young function. This phenomenon is in some sense magnified in the multi-phase framework: in \cite{deoh} it is observed that each $p_{\nu}$-phase interacts only with the elliptic $p$-phase as quantified by \eqref{bounds}; in particular no additional relation between $p_{\nu_{1}}$, $p_{\nu_{2}}$ or $\alpha_{\nu_{1}}$, $\alpha_{\nu_{2}}$ with $\nu_{1}\not =\nu_{2}\in \texttt{I}_{\kk}$ should be imposed. On a more formal level, according to the classification done in \cite{ciccio} we see that the integrand in $\mathcal{H}(\cdot)$ is \emph{pointwise uniformly elliptic}, in the sense that its ellipticity ratio is uniformly bounded:
\begin{eqnarray}\label{erh}
\mathcal{R}_{H}(z):=\sup_{x\in B}\frac{\mbox{highest eigenvalue of} \ \partial^{2}H(x,z)}{\mbox{lowest eigenvalue of} \ \partial^{2}H(x,z)}\le c(n,p,p_{1},\cdots,p_{\kk})
\end{eqnarray}
for any $z\in \mathbb{R}^{n}$ and all balls $B\Subset \Omega$. However, the possible vanishing of the coefficients creates a deficit in the structure that can be better measured via a nonlocal counterpart of the ellipticity ratio defined as
\begin{eqnarray*}
\mathfrak{R}_{H}(z):=\frac{\sup_{x\in B}\mbox{highest eigenvalue of} \ \partial^{2}H(x,z)}{\inf_{x\in B}\mbox{lowest eigenvalue of} \ \partial^{2}H(x,z)}\lesssim 1+\sum_{\nu=1}^{\kk}\nr{a_{\nu}}_{L^{\infty}(B)}\snr{z}^{p_{\nu}-p},
\end{eqnarray*}
which may blow up as $\snr{z}\to \infty$. From this analysis it is clear that nonuniform ellipticity of multi-phase integrands is caused by the coefficients, but it is rather soft and still allows a perturbative approach to regularity. The multi-phase energy is a particular instance of Musielak-Orlicz functional, an abstract class of variational integrals described for instance in \cite{hhb}, that permits to treat in a unified fashion the regularity of minima of several model functionals such as double phase, multi-phase, $p(x)$-Laplacian or double phase with variable exponent and the functional analytic properties of related Lagrangian spaces, see \cite{yags,ba2,bkhh,ck,cyz,HH,hhl,HO,HHO,KL} for an (incomplete) list of references and \cite{Majmaa,dark2} for reasonable surveys. It is worth mentioning that energy $\mathcal{H}(\cdot)$ also falls into the realm of functional with $(p,q)$-growth, i.e. variational integrals defined by means of a sufficiently smooth integrand $F\colon \Omega\times \mathbb{R}^{n}\to \mathbb{R}$ with a rate of nonuniform ellipticity stronger than \eqref{erh}, i.e: 
\begin{eqnarray*}
\begin{cases}
\ \snr{z}^{p}\lesssim F(x,z)\lesssim 1+\snr{z}^{q}\\
\ \mathcal{R}_{F}(z)\lesssim \snr{z}^{q-p}
\end{cases}\qquad  \mbox{with} \ \ 1<p\le q.
\end{eqnarray*}
This class of functionals has first been introduced in the seminal papers \cite{Mlav,M0,M1,M2} and intensively investigated since then in \cite{BM,BS,bc,bobr,BB,cdll,CKP,CKP1,choe,demii,ciccio,demib,eslemi,HS,S}, see also \cite{bi2,Majmaa,dark2} for an overview of the state of the art. The main idea in this case consists in neglecting the precise structure of the integrand and retaining only the extremals of the growth. In such way it is possible to prove regularity results for minima of quite a large family of variational integrals at the price of imposing precise closeness conditions between exponents $(p,q)$ and loosing some informations that are distinctive of the specific structure, compare in this perspective \cite[Theorem 1]{ciccio} with \cite[Theorem 3]{ciccio}. The regularity for general functionals with $(p,q)$-growth is guaranteed provided that $q/p\le 1+\texttt{o}(n)$, where $\texttt{o}(n)\to_{n\to \infty}0$. This turns out to be a necessary and sufficient condition for regularity, cf. \cite{eslemi,M1,b32} about the counterexamples/sharpness of the upper bound on the ratio $q/p$ and \cite{BS,HS,S} for improvements in the autonomous setting. The constraint linking exponents $(p,q)$ has interpolative nature in the sense that if minimizers a priori feature a higher regularity than the one naturally allowed by the ellipticity of the functional, then the restriction imposed on the size of $q/p$ can be relaxed, in particular it can be made independent on the space dimension, cf. \cite{bacomi,bi2,cdll,CKP,choe,comib,demii,demima,demib,ok}. The main tool exploited in most of such papers are Gagliardo-Nirenberg type inequalities \cite{BrMi} that grant a trading between the extra regularity properties of minima and the higher integrability of their gradients. This transaction weakens in some sense the nonuniform ellipticity of the functional, thus either allowing for larger bounds on $q/p$ or drastically reducing the rate of fractional differentiability of the gradient needed for boost its integrability. The latter is the cornerstone of the arguments presented here.
\subsubsection*{Organization of the paper} This paper is organized as follows. In Section \ref{pre} we describe our notation and collect some auxiliary results, Section \ref{rmp} contains an overview of the regularity theory for local minimizers of multi-phase integrals and Sections \ref{crh}-\ref{czss} are devoted to the proofs of Theorems \ref{revh}-\ref{czt} respectively.
\section{Preliminaries}\label{pre}
In this section we shall collect some well-known results that will be useful in the proof of Theorems \ref{revh}-\ref{czt}.
\subsection{Notation}\label{not}
We denote by $\Omega\subset \er^n$ an open domain and, since our estimates will be local, we shall always assume, without loss of generality, that $\Omega$ is also bounded. We denote by $c$ a general constant larger than one. Different occurrences from line to line will be still denoted by $c$. Important dependencies on parameters will be as usual emphasized by putting them in parentheses. We shall denote $\N$ as the set of positive integers. As usual, we denote by $ B_r(x_0):= \{x \in \er^n  :   |x-x_0|< r\}$ the open ball with center $x_0$ and radius $r>0$; when it is clear from the context, we omit denoting the center, i.e., $B_r \equiv B_r(x_0)$. When not otherwise stated, different balls in the same context will share the same center. Finally, with $B$ being a given ball with radius $r$ and $\gamma$ being a positive number, we denote by $\gamma B$ the concentric ball with radius $\gamma r$ and by $B/\gamma \equiv (1/\gamma)B$. In denoting several function spaces like $L^p(\Omega), W^{1,p}(\Omega)$, we shall denote the vector valued version by $L^p(\Omega,\er^k), W^{1,p}(\Omega,\er^k)$ in the case the maps considered take values in $\er^k$, $k\in \en$. With $\mathcal B \subset \er^{n}$ being a measurable subset with bounded positive measure $0<|\mathcal B|<\infty$, and with $g \colon \mathcal B \to \er^{k}$, $k\geq 1$, being a measurable map, we shall denote the integral average of $g$ over $\mathcal B$ by  
$$
   (g)_{\mathcal B} \equiv \mint_{\mathcal B}  g(x) \ \dx  :=  \frac{1}{|\mathcal B|}\int_{\mathcal B}  g(x) \ \dx\;.
$$
Moreover, if $g\colon \Omega\to \mathbb{R}^{k}$ is any map, $U\subset \Omega$ is an open set and $\beta \in (0,1]$ is a given number we shall denote
\begin{flalign*}
[g]_{0,\beta;U}:=\sup_{x,y \in U; x\not=y}\frac{\snr{g(x)-g(y)}}{\snr{x-y}^{\beta}}, \qquad [g]_{0,\beta}:=[g]_{0,\beta;\Omega}.
\end{flalign*}
The quantity in the previous definition is a seminorm and $g$ is included in the H\"older space $C^{0,\beta}(U,\mathbb{R}^{k})$ iff $[g]_{0,\beta;U}<\infty$. We also point out that $g\in C^{1,\beta}(U,\mathbb{R}^{k})$ provided that $Dg\in C^{0,\beta}(U,\mathbb{R}^{k\times n})$. Finally, for the sake of simplicity, we collect the main parameters of the problem in the shorthands
\begin{eqnarray*}
\begin{cases}
\ \texttt{A}:=\max_{\nu\in \texttt{I}_{\kk}}\nr{a_{\nu}}_{L^{\infty}(\Omega)}\\
\ \texttt{data}_{0}:=(n,p,p_{1},\cdots,p_{\kk},\alpha_{1},\cdots,\alpha_{\kk}),\\
\ \texttt{data}:=(\texttt{data}_{0},[a_{1}]_{0,\alpha_{1}},\cdots, [a_{\kk}]_{0,\alpha_{\kk}},\kk)\\
\ \texttt{data}_{\textnormal{cz}}:=(\texttt{data},\texttt{A},\Lambda,\nr{H(\cdot,Du)}_{L^{1}(\ti{\Omega}_{0})},\nr{H(\cdot,F)}_{L^{\gamma}(\ti{\Omega}_{0})},\gamma,\dist(\ti{\Omega}_{0},\partial \Omega)),
\end{cases}
\end{eqnarray*}
see Sections \ref{rmp}-\ref{czss} for more informations about the quantities mentioned in the previous display.
\subsection{On fractional Sobolev spaces}
Given a function $w \colon \Omega \to \mathbb{R}^{k}$, $k\ge 1$ and a vector $h \in \mathbb{R}^n$, we denote by $\tau_{h}\colon L^1(\Omega,\mathbb{\er}^{k}) \to L^{1}(\Omega_{|h|},\mathbb{R}^{k})$ the standard finite difference operator pointwise defined as
\begin{flalign*}
\tau_{h}w(x):=w(x+h)-w(x) \quad \mbox{for a.e.} \ x\in \Omega_{\snr{h}},
\end{flalign*}
where $\Omega_{|h|}:=\{x \in \Omega \, : \, 
\dist(x, \partial \Omega) > |h|\}$. Let us record the fundamentals of fractional Sobolev spaces, see \cite{pala} for more details on this matter.
\begin{definition}\label{fra1def}
Let $\Omega\subset \mathbb{R}^{n}$ be an open set with $n\ge 2$ (the case $\Omega\equiv \mathbb{R}^{n}$ is allowed as well), $\alpha\in (0,1)$, $p\in [1,\infty)$ and $k\in \N$ be numbers. The fractional Sobolev space $W^{\alpha ,p}(\Omega,\mathbb{R}^k )$ is defined prescribing that $w \colon \Omega \to \mathbb{R}^k$ belongs to  $W^{\alpha ,p}(\Omega,\mathbb{R}^{k} )$ iff the following Gagliardo type norm is finite:
\begin{flalign*}
\nr{w}_{W^{\alpha,p}(\Omega)}&:=\nr{w}_{L^{p}(\Omega)}+\left(\int_{\Omega} \int_{\Omega}  
\frac{\snr{w(x)
- w(y) }^{p}}{\snr{x-y}^{n+\alpha p}} \ \dx \, \dy \right)^{1/p}=:\nr{w}_{L^{p}(\Omega)}+[w]_{\alpha,p;\Omega}.
\end{flalign*}
Accordingly, if $\alpha = [\alpha]+\{\alpha\}\in \en + (0,1)>1$, we say that $w\in W^{\alpha ,p}(\Omega,\mathbb{R}^k )$ iff the following quantity is finite
\begin{flalign*}
\nr{w}_{W^{\alpha,p}(\Omega)}:=\nr{w}_{W^{[\alpha],p}(\Omega)}+[D^{[\alpha]}w]_{\{\alpha\},p;\Omega}.
\end{flalign*}
The local variant $W^{\alpha ,p}_{loc}(\Omega,\er^k )$ is defined by requiring that $w \in W^{\alpha ,p}(\tilde{\Omega},\mathbb{R}^k)$ for every open subset $\tilde{\Omega} \Subset \Omega$. 
\end{definition}
A class of spaces that is strictly related to fractional Sobolev spaces is that of Nikol'skii spaces.
\begin{definition}\label{fra2def}
Let $\Omega\subset \mathbb{R}^{n}$ be an open set with $n\ge 2$ and $\alpha \in (0,1)$, $p\in [1, \infty)$, $k \in \mathbb{R}^{n}$ be numbers. The Nikol'skii space $N^{\alpha,p}(\Omega,\mathbb{R}^k )$ is defined prescribing that $w\in N^{\alpha,p}(\Omega,\mathbb{R}^k )$ iff 
$$\nr{w}_{N^{\alpha,p}(\Omega )} :=\nr{w}_{L^{p}(\Omega)} + \left(\sup_{\snr{h}\not=0}\, \int_{\Omega_{|h|}} 
\frac{\snr{w(x+h)
- w(x) }^{p}}{|h|^{\alpha p}} \ dx  \right)^{1/p}\;.$$
The local variant $N^{\alpha,p}_{\loc}(\Omega,\mathbb{R}^{k} )$ is defined by requiring that $w \in N^{\alpha,p}(\tilde \Omega,\mathbb{R}^{k})$ for every open subset $\tilde{\Omega} \Subset \Omega$.
\end{definition}
Whenever $\Omega$ is a sufficiently regular domain, it is $W^{\alpha_{0},p}(\Omega,\mathbb{R}^{k})\nsubseteq N^{\alpha_{0},p}(\Omega,\mathbb{R}^{k})\nsubseteq W^{\beta,p}(\Omega,\mathbb{R}^{k})$ for all $\beta\in (0,\alpha_{0})$. This chain of inclusions can be in some sense quantified, and this is the content of the next lemma, cf. \cite[Section 2.2]{demii}.
\begin{lemma}\label{l2}
Let $B_{r}\Subset \er^n$ be a ball with $r\leq 1$, $w\in L^{p}(B_{r},\mathbb{R}^{k})$, $p>1$ and assume that, for $\alpha \in (0,1]$, $S\ge 1$ and concentric balls $B_{\rr}\Subset B_{r}$, there holds
 $$
\nr{\tau_{h}w}_{L^{p}(B_{\rr},\er^k)}\le S\snr{h}^{\alpha } \quad \mbox{
for every $h\in \mathbb{R}^{n}$ with $0<\snr{h}\le \frac{r-\rr}{K}$, where $K \geq 1$}\;.$$ 
Then $w\in W^{\beta,p}(B_{\rr},\mathbb{R}^{k})$ whenever $\beta\in (0,\alpha )$ and
$$
\nr{w}_{W^{\beta,p}(B_{\rr},\er^k)}\le\frac{c}{(\alpha -\beta)^{1/p}}
\left(\frac{r-\rr}{K}\right)^{\alpha -\beta}S+c\left(\frac{K}{r-\rr}\right)^{n/p+\beta} \nr{w}_{L^{p}(B_{r},\er^k)}\;,
$$
holds, where $c\equiv c(n,p)$. 
\end{lemma}
We conclude this section with a fractional Gagliardo-Nirenberg type inequality, whose proof can be found in \cite[Corollary 3.2]{BrMi}, see also \cite[Lemma 2.6]{comib} for a localized version. 
\begin{lemma}\label{fraclem}
Let $B_{\rr}\Subset B_{r}\Subset \er^n$ be concentric balls with $r \leq 1$. Let $0< s_{1}<1<s_{2}<2$, $1<p,q<\infty$, $t>1$ and $\theta\in (0,1)$ be such that
\begin{flalign*}
1=\theta s_{1}+(1-\theta)s_{2},\qquad \frac{1}{t}=\frac{\theta}{q}+\frac{1-\theta}{p}.
\end{flalign*}
Then every function $w\in W^{s_{1},q}(B_{r})\cap W^{s_{2},p}(B_{r})$ belongs to $W^{1,t}(B_{\rr})$ and the inequality 
\eqn{viavia}
$$
\nr{Dw}_{L^{t}(B_{\rr})}\le \frac{c}{(r-\rr)^{\kappa}}[w]^{\theta}_{s_{1},q;B_{r}}\nr{Dw}^{1-\theta}_{W^{s_{2}-1,p}(B_{r})}
$$
holds for constants $c,\kappa\equiv c,\kappa(n,s_{1},s_{2},p,q,t)$.
\end{lemma}
\subsection{Tools for nonlinear problems}
When dealing with $m$-Laplacean type problems with $m>1$, we shall often use the auxiliary vector fields $V_{m}\colon \er^{n} \to  \er^{n}$, defined by
\begin{flalign*}
V_{m}(z):= |z|^{(p-2)/2}z, \qquad m\in (1,\infty)
\end{flalign*}
whenever $z \in \er^{n}$. In Sections \ref{crh}-\ref{czss}, we shall adopt the above definition with $m\in \{p,p_{1},\cdots,p_{\kk}\}$. A useful related inequality is contained in the following
\begin{flalign}\label{vpvp}
\snr{V_{m}(z_{1})-V_{m}(z_{2})}\approx (\snr{z_{1}}^{2}+\snr{z_{2}}^{2})^{(p-2)/4}\snr{z_{1}-z_{2}}, 
\end{flalign}
where the equivalence holds up to constants depending only on $n,m$. An important property which is usually related to such field is recorded in the following lemma.
\begin{lemma}\label{l6}
Let $t>-1$, and $z_{1},z_{2}\in \mathbb{R}^{n}$ be so that $\snr{z_{1}}+\snr{z}_{2}>0$. Then
\begin{flalign*}
\int_{0}^{1}\snr{z_{1}+\lambda(z_{2}-z_{1})}^{t} \ \d\lambda\sim (\snr{z_{1}}^{2}+\snr{z_{2}}^{2})^{\frac{t}{2}},
\end{flalign*}
with constants implicit in "$\sim$" depending only from $t$.
\end{lemma}
Finally, the "simple, but fundamental" iteration lemma of \cite[Section 1]{gigi}. 
\begin{lemma}\label{iter}
Let $\mathcal{Z}\colon [\rr,R)\to [0,\infty)$ be a function which is bounded on every interval $[\varrho, R_*]$ with $R_*<R$. Let $\varepsilon\in (0,1)$, $a_1,a_2,\gamma_{1},\gamma_{2}\ge 0$ be numbers. If
\begin{flalign*}
\mathcal{Z}(\tau_1)\le \varepsilon \mathcal{Z}(\tau_2)+ \frac{a_1}{(\tau_2-\tau_1)^{\gamma_{1}}}+\frac{a_2}{(\tau_2-\tau_1)^{\gamma_{2}}}\ \ \mbox{for all} \ \rr\le \tau_1<\tau_2< R\;,
\end{flalign*}
then
\begin{flalign*}
\mathcal{Z}(\rr)\le c\left[\frac{a_1}{(R-\rr)^{\gamma_{1}}}+\frac{a_2}{(R-\rr)^{\gamma_{2}}}\right]\;,
\end{flalign*}
holds with $c\equiv c(\varepsilon,\gamma_{1},\gamma_{2})$.
\end{lemma}

\section{Regularity theory for local minimizers of multi-phase problems}\label{rmp}
In this section we collect some well-known regularity results for minima of functional $\mathcal{H}(\cdot)$, i.e. maps verifying the following definition.
\begin{definition}\label{min}
A function $v\in W^{1,p}_{\loc}(\Omega)$ is a local minimizer of functional $\mathcal{H}(\cdot)$ if for every open subset $\tilde{\Omega}\Subset \Omega$ we have $\mathcal{H}(v,\ti{\Omega})<\infty$ and $\mathcal{H}(v,\ti{\Omega})\le \mathcal{H}(w,\ti{\Omega})$ for any competitor $w\in v+W^{1,p}_{0}(\ti{\Omega})$ so that $H(\cdot,Dw)\in L^{1}(\ti{\Omega})$.
\end{definition}
The details of the proof of all the results listed below can be found in \cite{deoh} for the case of three phases, i.e. $H(x,z)\equiv \left[\snr{z}^{p}+a_{1}(x)\snr{z}^{p_{1}}+a_{2}(x)\snr{z}^{p_{2}}\right]$, but, as stressed in \cite[Section 1]{deoh}, they can be adapted in a straightforward way to an arbitrary number of phases, see also \cite[Section 2 and Theorems 7.2-7.4]{HO}. We start by discussing a peculiar feature of variational integrals with Musielak-Orlicz structure which is the absence of Lavrentiev Phenomenon, see \cite[Lemma 1]{demicz} and \cite[Lemma 13]{eslemi}.
\begin{lemma}\label{lgh}
Under assumptions \eqref{ab}-\eqref{bounds}, let $w\in W^{1,p}_{\loc}(\Omega)$ be any function so that whenever $B\Subset \Omega$ is a bounded, open set it is $\nr{H(\cdot,Dw)}_{L^{1+\delta'}(B)}<\infty$ for some $\delta'>0$. Then there exists a sequence of smooth maps $\{\ti{w}_{j}\}_{j\in \N}\subset C^{\infty}_{\loc}(\Omega)$ so that it holds
\begin{flalign}\label{0}
\begin{cases}
\ \ti{w}_{j}\to w \ \ \mbox{in} \ \ W^{1,p(1+\delta')}(B)\\
\ \nr{H(\cdot,D\ti{w}_{j})}_{L^{1}(B)}\to \nr{H(\cdot,Dw)}_{L^{1}(B)}\\
\ \nr{H(\cdot,D\ti{w}_{j})}_{L^{1+\delta'}(B)}\to \nr{H(\cdot,w)}_{L^{1+\delta'}(B)}.
\end{cases}
\end{flalign}
\end{lemma}
Next, a Sobolev-Poincar\'e inequality for multi-phase problems, \cite[Lemma 2]{deoh}.
\begin{lemma}
Under assumptions \eqref{ab}-\eqref{bounds}, let $B_{\rr}\Subset \mathbb{R}^{n}$ be a ball with radius $\rr\in (0,1]$ and $w\in W^{1,p}(B_{\rr})$ be any function so that $H(\cdot,Dw)\in L^{1}(B_{\rr})$. Then there are a positive constant $c\equiv c(\texttt{data}_{0})$ and an exponent $d\equiv d(n,p,p_{1},\cdots,p_{\kk})\in (0,1)$ so that
\begin{flalign}\label{sopo}
\mint_{B_{\rr}}H\left( x, \frac{w-(w)_{B_{\rr}}}{\rr} \right) \ \dx\le c\left(1+\sum_{\nu=1}^{\kk}[a_{\nu}]_{0,\alpha_{\nu};B_{\rr}}\nr{Dw}_{L^{p}(B_{\rr})}^{p_{\nu}-p}\right)\left(\mint_{B_{\rr}}H(x,Dw)^{d} \ \dx\right)^{1/d}.
\end{flalign}
\end{lemma}
Let us record a local higher integrability result of Gehring type, cf. \cite[Lemma 4]{deoh}.
\begin{lemma}\label{lgeh}
Under assumptions \eqref{ab}-\eqref{bounds}, let $B_{\sigma}\Subset \Omega$ be any ball with radius $\sigma\in (0,1]$ and $v\in W^{1,p}_{\loc}(\Omega)$ be a local minimizer of functional $\mathcal{H}(\cdot)$ so that $\nr{H(\cdot,Dv)}_{L^{1}(B_{\sigma})}\le M$ for some constant $M>0$. Then there exists a positive higher integrability threshold $\delta_{g}\equiv \delta_{g}(\texttt{data},M)$ so that
\begin{eqnarray}\label{geh}
\left(\mint_{B_{\sigma/2}}H(x,Dv)^{1+\delta} \ \dx\right)^{\frac{1}{1+\delta}}\le c\mint_{B_{\sigma}}H(x,Dv) \ \dx,
\end{eqnarray}
for all $\delta\in (0,\delta_{g}]$, with $c\equiv c(\texttt{data},M)$. 
\end{lemma}
The global counterpart of Lemma \ref{lgh} is in the next lemma.
\begin{lemma}\label{bgeh}
Under assumptions \eqref{ab}-\eqref{bounds}, let $B\Subset \Omega$ be a ball with radius $\texttt{r}(B)\in (0,1]$, $u_{0}\in W^{1,p(1+\delta_{0})}(B)$ for some $\delta_{0}>0$ with $\nr{H(\cdot,Du_{0})}_{L^{1}(B)}\le M_{0}$ be any function and $v_{0}\in u_{0}+W^{1,p}_{0}(B)$ be the solution of Dirichlet problem
\begin{eqnarray*}
u_{0}+W^{1,p}_{0}(B)\ni w\mapsto \min \mathcal{H}(w,B).
\end{eqnarray*}
There exists an higher integrability threshold $\sigma_{g}\equiv \sigma_{g}(\texttt{data},M_{0},\delta_{0})\in (0,\delta_{0})$ so that
\begin{eqnarray*}
\mint_{B}H(x,Dv_{0})^{1+\sigma_{g}} \ \dx\le c\mint_{B}H(x,Du_{0})^{1+\sigma_{g}} \ \dx,
\end{eqnarray*}
for $c\equiv c(\texttt{data},M_{0},\delta_{0})$.
\end{lemma}
We further recall a straightforward manipulation of \cite[Theorem 2]{deoh}.
\begin{theorem}\label{mor}
Under assumptions \eqref{ab}-\eqref{bounds}, let $B\Subset \Omega$ be an open, bounded set and $v\in W^{1,p}_{\loc}(\Omega)$ be a local minimizer of functional $\mathcal{H}(\cdot)$ so that $\nr{H(\cdot,Dv)}_{L^{1+\delta_{g}}(B)}\le M_{g}$, where $\delta_{g}$ is the higher integrability threshold coming from Lemma \ref{lgeh}. Then, whenever $B_{\sigma_{1}}\subset B_{\sigma_{2}}\Subset B$ are concentric balls with radii $0<\sigma_{1}\le \sigma_{2}\le 1$, for every $\beta\in (0,n)$ it holds that
\begin{eqnarray*}
\int_{B_{\sigma_{1}}}H(x,Dv) \ \dx \le c\left(\frac{\sigma_{1}}{\sigma_{2}}\right)^{n-\beta}\int_{B_{\sigma_{2}}}H(x,Dv) \ \dx,
\end{eqnarray*}
with $c\equiv c(\texttt{data},M_{g},\beta)$. In particular, 
$v\in C^{0,\gamma_{0}}(B)$ for all $\gamma_{0}\in (0,1)$ with
\begin{eqnarray}\label{hh}
[v]_{0,\gamma_{0};B_{2\sigma/3}}\le c \sigma^{1-\gamma_{0}}\left(\mint_{B_{\sigma}}H(x,Dv) \ \dx\right)^{1/p},
\end{eqnarray}
for $c\equiv c(\texttt{data},M_{g},\gamma_{0})$.
\end{theorem}
Finally, we conclude this section with the main result of \cite{deoh}.
\begin{theorem}\label{ht}
Let $v\in W^{1,p}_{\loc}(\Omega)$ be a local minimizer of functional $\mathcal{H}(\cdot)$, with \eqref{ab}-\eqref{bounds} in force. Then $v\in C^{1,\beta_{0}}_{\loc}(\Omega)$ for some $\beta_{0}\equiv \beta_{0}(\texttt{data}_{0})$.
\end{theorem}

\begin{remark}\label{inc}
\emph{We stress that all the constants appearing in Lemmas \ref{lgeh}-\ref{bgeh} and Theorem \ref{mor} are nondecreasing in $M$, $M_{0}$ and $M_{g}$ respectively, cf. \cite{deoh}.}
\end{remark}

\section{Conditional reverse H\"older inequalities}\label{crh}
In this section we prove our main result, i.e. a reverse H\"older inequality for minima of $\mathcal{H}(\cdot)$ in the spirit of those appearing in \cite{comib,comi} without imposing any relation between the H\"older continuity exponents $\{\alpha_{\nu}\}_{\nu=1}^{\kk}$. A similar result has been obtained in \cite[Theorem 4.1]{bbo} for generalized multi-phase problems with the additional technical constraint \eqref{ns}. We believe that our proof can be adapted to more general functionals than $\mathcal{H}(\cdot)$ that still preserve specific Musielak-Orlicz structure. Following a by now standard terminology see \cite{bacomi,comib,comicz,comi,demicz,ciccio,demima} and in particular \cite[Section 4]{deoh}, given any ball $B_{\rr}(x_{0})\Subset \Omega$, we identify three scenarios, according to the behavior of coefficients $\{a_{\nu}(\cdot)\}_{\nu=1}^{\kk}$. Precisely, given any constant $J\ge 4$, we shall say that $\mathcal{H}(\cdot)$ is in the \emph{degenerate} phase $\texttt{deg}_{\textnormal{J}}(B_{\rr}(x_{0}))$ on $B_{\rr}(x_{0})$ if
\begin{eqnarray*}
\inf_{x\in B_{\rr}(x_{0})}a_{\nu}(x)\le J[a_{\nu}]_{0,\alpha_{\nu};B_{\rr}(x_{0})}\rr^{\alpha_{\nu}}\qquad \mbox{for all} \ \ \nu\in \texttt{I}_{\kk},
\end{eqnarray*}
while $\mathcal{H}(\cdot)$ is in the \emph{nondegenerate} phase $\texttt{ndeg}_{\textnormal{J}}(B_{\rr}(x_{0}))$ when
\begin{eqnarray*}
\inf_{x\in B_{\rr}(x_{0})}a_{\nu}(x)> J[a_{\nu}]_{0,\alpha_{\nu};B_{\rr}(x_{0})}\rr^{\alpha_{\nu}}\qquad \mbox{for all} \ \ \nu\in \texttt{I}_{\kk},
\end{eqnarray*}
while $\mathcal{H}(\cdot)$ is in a \emph{mixed} phase $\texttt{mix}_{\textnormal{J}}(B_{\rr}(x_{0}))$ provided that the set of indexes $\texttt{I}_{\kk}$ is the union of two nonempty subsets $\texttt{d},\texttt{nd}\subset \texttt{I}_{\kk}$ which can be characterized as
\begin{eqnarray*}
\begin{cases}
\ \inf_{x\in B_{\rr}(x_{0})}a_{\nu}(x)\le J[a_{\nu}]_{0,\alpha_{\nu};B_{\rr}(x_{0})}\rr^{\alpha_{\nu}}\qquad \mbox{for all} \ \ \nu\in\texttt{d}\\
\ \inf_{x\in B_{\rr}(x_{0})}a_{\nu}(x)> J[a_{\nu}]_{0,\alpha_{\nu};B_{\rr}(x_{0})}\rr^{\alpha_{\nu}}\qquad \mbox{for all} \ \ \nu\in\texttt{nd}.
\end{cases}
\end{eqnarray*}
The above configurations will play a key role in the next sections.
\subsection{Proof of Theorem \ref{revh}}
For the transparency of presentation, we split the proof of Theorem \ref{revh} into nine steps. Since the dependencies of the constants declared throughout the proof may seem quite weird, we shall provide a detailed explanation of the behavior of such constants in \emph{Step 9}.
\subsubsection*{Step 1: scaling and approximation}
Let $v\in W^{1,p}_{\loc}(\Omega)$ be a local minimizer of functional $\mathcal{H}(\cdot)$ and $B_{\rr}(x_{0})\Subset \Omega$ be any ball with radius $\rr\in (0,1]$. By Lemma \ref{lgeh} we know that $H(\cdot,Dv)\in L^{1+\delta_{g}}(B_{\rr}(x_{0}))$ for some $\delta_{g}\equiv \delta_{g}(\texttt{data},\nr{H(\cdot,Dv)}_{L^{1}(B_{\rr}(x_{0}))})$, so Lemma \ref{lgh} applies and we obtain a sequence $\{\ti{v}_{j}\}_{j\in \N}\subset C^{\infty}(\bar{B}_{\rr}(x_{0}))$ so that \eqref{0} holds with $B\equiv B_{\rr}(x_{0})$. We blow $v$ on $B_{\rr}(x_{0})$ by defining $B_{1}(0)\ni x\mapsto v_{\rr}(x):=(v(x_{0}+\rr x)-(v)_{B_{\rr}(x_{0})})\rr^{-1}$ and notice that a simple scaling argument shows that $v_{\rr}\in W^{1,p}(B_{1}(0))$ is a local minimizer of functional
\begin{eqnarray*}
W^{1,p}(B_{1}(0))\ni w\mapsto\mathcal{H}_{\rr}(w,B_{1}(0)):=\int_{B_{1}(0)}H_{\rr}(x,z) \ \dx,
\end{eqnarray*}
with $B_{1}(0)\ni x\mapsto a_{\nu,\rr}(x):=a_{\nu}(x_{0}+\rr x)$ for all $\nu\in \texttt{I}_{\kk}$ and
\begin{eqnarray*}
 H_{\rr}(x,z):=\left[\snr{z}^{p}+\sum_{\nu=1}^{\kk}a_{\nu,\rr}(x)\snr{z}^{p_{\nu}}\right].
\end{eqnarray*}
By definition we have that
\begin{eqnarray}\label{0.1}
\begin{cases}
\ \nr{a_{\nu,\rr}}_{L^{\infty}(B_{1}(0))}=\nr{a_{\nu}}_{L^{\infty}(B_{\rr}(x_{0}))}\quad &\mbox{for all} \ \ \nu\in\texttt{I}_{\kk}\\
\ [a_{\nu,\rr}]_{0,\alpha_{\nu};B_{1}(0)}=\rr^{\alpha_{\nu}}[a_{\nu}]_{0,\alpha_{\nu};B_{\rr}(x_{0})}\quad &\mbox{for all} \ \ \nu\in \texttt{I}_{\kk}\\
\ \mathcal{H}_{\rr}(v_{\rr},B_{1}(0))=\rr^{-n}\mathcal{H}(v,B_{\rr}(x_{0})).
\end{cases}
\end{eqnarray}
We stress that by construction, $v_{\rr}$ retains the same higher integrability features of $v$, i.e. $H_{\rr}(\cdot,Dv_{\rr})\in L^{1+\delta_{g}}(B_{1}(0))$ where $\delta_{g}\equiv \delta_{g}(\texttt{data},\nr{H(\cdot,Dv)}_{L^{1}(B_{\rr}(x_{0}))})$ is the same higher integrability exponent of $v$. Moreover, setting $B_{1}(0)\ni x\mapsto\ti{v}_{j,\rr}(x):=(\ti{v}_{j}(x_{0}+\rr x)-(\ti{v}_{j})_{B_{\rr}(x_{0})})\rr^{-1}$, by \eqref{0} with $B\equiv B_{\rr}(x_{0})$ we have a sequence $\{v_{j,\rr}\}_{j\in \N}\subset C^{\infty}(\bar{B}_{1}(0))$ so that
\begin{flalign}\label{0r}
\begin{cases}
\ \ti{v}_{j,\rr}\to v_{\rr} \ \ \mbox{strongly in} \ \ W^{1,p(1+\delta_{g})}(B_{1}(0))\\
\ \nr{H_{\rr}(\cdot,D\ti{v}_{j,\rr})}_{L^{1}(B_{1}(0))}\to \nr{H_{\rr}(\cdot,Dv_{\rr})}_{L^{1}(B_{1}(0))}\\
\ \nr{H_{\rr}(\cdot,D\ti{v}_{j,\rr})}_{L^{1+\delta_{g}}(B_{1}(0))}\to \nr{H_{\rr}(\cdot,Dv_{\rr})}_{L^{1+\delta_{g}}(B_{1}(0))}.
\end{cases}
\end{flalign}
For $\nu\in \texttt{I}_{\kk}$ and $j\in \N$, we correct the growth of $H_{\rr}(\cdot)$ by introducing the regularized integrands
\begin{eqnarray*}
H_{j}(x,z):=H_{\rr}(x,z)+\sum_{\nu=1}^{\kk}\sigma_{j}^{\nu}\snr{z}^{p_{\nu}}\equiv \snr{z}^{p}+\sum_{\nu=1}^{\kk}\left(a_{\nu,\rr}(x)+\sigma^{\nu}_{j}\right)\snr{z}^{p_{\nu}}, 
\end{eqnarray*}
where we set
$$
\sigma_{j}^{\nu}:=j^{-1}\left(1+j+\nr{D\ti{v}_{j,\rr}}_{L^{p_{\nu}}(B_{1}(0))}^{2p_{\nu}}+\nr{D\ti{v}_{j,\rr}}_{L^{p_{\nu}(1+\delta_{g})}(B_{1}(0))}^{2p_{\nu}}\right)^{-1}.
$$
By very definition, it is
\begin{eqnarray}\label{1}
\sum_{\nu=1}^{\kk}\sigma^{\nu}_{j}\int_{B_{1}(0)}\snr{D\tilde{v}_{j,\rr}}^{p_{\nu}} \ \dx+\sum_{\nu=1}^{\kk}(\sigma^{\nu}_{j})^{1+\delta_{g}}\int_{B_{1}(0)}\snr{D\ti{v}_{j,\rr}}^{p_{\nu}(1+\delta_{g})} \ \dx\to 0
\end{eqnarray}
Keeping in mind $\eqref{ab}_{2}$, we set $\bar{p}:=\max_{\nu\in \texttt{I}_{\kk}}p_{\nu}$ and define the family of auxiliary Dirichlet problems 
\begin{eqnarray}\label{pda}
\ti{v}_{j,\rr}+W^{1,\bar{p}}_{0}(B)\ni w\mapsto \mathcal{H}_{j}(w,B_{1}(0)):=\int_{B_{1}(0)}H_{j}(x,z) \ \dx.
\end{eqnarray}
Direct methods assure that problem \eqref{pda} admits a unique solution $v_{j}\in \ti{v}_{j,\rr}+W^{1,p}_{0}(B_{1}(0))$ and, according to the regularity theory in \cite{li1} it is
\begin{eqnarray}\label{d3}
v_{j}\in W^{1,\infty}(B_{1}(0)),
\end{eqnarray}
given that $\ti{v}_{j,\rr}\in C^{\infty}(\bar{B}_{1}(0))$ and $\sigma^{\nu}_{j}>0$ for all $\nu\in \texttt{I}_{\kk}$, so $H_{j}(\cdot)$ has standard $\bar{p}$-growth. We further notice that functional $\mathcal{H}_{j}(\cdot)$ is of multi-phase type. In fact \eqref{bounds} is always in force and \eqref{ab} trivially holds since the coefficients $a_{\nu,\rr}+\sigma^{\nu}_{\rr}\in C^{0,\alpha_{\nu}}(\Omega)$ verify $[a_{\nu,\rr}+\sigma^{\nu}_{\rr}]_{0,\alpha_{\nu};B_{1}(0)}=[a_{\nu,\rr}]_{0,\alpha_{\nu};B_{1}(0)}$ for all $\nu\in \texttt{I}_{\kk}$, therefore Lemma \ref{bgeh} applies and there is an exponent $\sigma_{g}\equiv \sigma_{g}(\texttt{data},\nr{H(\cdot,Dv)}_{L^{1}(B_{\rr}(x_{0}))})\in (0,\delta_{g})$ so that
\begin{eqnarray}
\nr{H_{j}(\cdot,Dv_{j})}_{L^{1+\sigma_{g}}(B_{1}(0))}&\le& c\nr{H_{j}(\cdot,D\ti{v}_{j,\rr})}_{L^{1+\sigma_{g}}(B_{1}(0))}\nonumber \\
&\stackrel{\eqref{0r},\eqref{1}}{\le}&c\left(\nr{H_{\rr}(\cdot,Dv_{\rr})}_{L^{1+\delta_{g}}(B_{1}(0))}+1\right),\label{0.2}
\end{eqnarray}
with $c\equiv c(\texttt{data},\nr{H(\cdot,Dv)}_{L^{1}(B_{\rr}(x_{0}))})$.
\subsubsection*{Step 2: covergence} Let us prove that the sequence $\{v_{j}\}_{j\in \N}\subset W^{1,\bar{p}}(B_{1}(0))\cap W^{1,\infty}(B_{1}(0))$ of solutions to problem \eqref{pda} converge to $v_{\rr}$, local minimizer on $B_{1}(0)$ of $\mathcal{H}_{\rr}(\cdot)$. By minimality it is
\begin{eqnarray}\label{2}
\mathcal{H}_{j}(v_{j},B_{1}(0))&\le& \mathcal{H}_{j}(\ti{v}_{j,\rr},B_{1}(0))\nonumber \\
&\le& \mathcal{H}_{\rr}(\ti{v}_{j,\rr},B_{1}(0))+\sum_{\nu=1}^{\kk}\int_{B_{1}(0)}\sigma^{\nu}_{j}\snr{D\ti{v}_{j,\rr}}^{p_{\nu}} \ \dx\nonumber \\
&\stackrel{\eqref{1}}{\le}&\mathcal{H}_{\rr}(\ti{v}_{j,\rr},B_{1})+\texttt{o}(j)\stackrel{\eqref{0r}_{2}}{\le}\mathcal{H}_{\rr}(v_{\rr},B_{1}(0))+\texttt{o}(j),
\end{eqnarray}
which means that (keep \eqref{0r}$_{1}$ in mind)
\begin{eqnarray}\label{3}
v_{j}\rightharpoonup \ti{v} \ \ \mbox{weakly in} \ \ W^{1,p}(B_{1}(0))\qquad \mbox{and}\qquad \left. \ti{v}\right|_{\partial B_{1}(0)}=\left. v\right|_{\partial B_{1}(0)}.
\end{eqnarray}
The content of the previous display allows using weak lower semicontinuity in \eqref{2} to get
\begin{eqnarray*}
\mathcal{H_{\rr}}(v_{\rr},B_{1}(0))&\stackrel{\eqref{3}_{2}}{\le}&\mathcal{H}_{\rr}(\ti{v},B_{1}(0))\le\liminf_{j\to \infty}\mathcal{H}_{\rr}(v_{j},B_{1}(0))\nonumber \\
&\le&\limsup_{j\to \infty}\mathcal{H}_{\rr}(v_{j},B_{1}(0)) \le\limsup_{j\to \infty}\mathcal{H}_{j}(v_{j},B_{1}(0))\nonumber \\
&\le& \limsup_{j\to \infty}\left[\mathcal{H}_{\rr}(\ti{v}_{j,\rr},B_{1}(0))+\sum_{\nu=1}^{\kk}\int_{B_{1}(0)}\sigma^{\nu}_{j}\snr{D\ti{v}_{j,\rr}}^{p_{\nu}} \ \dx\right]\nonumber \\
&\stackrel{\eqref{1}}{\le}&\mathcal{H}_{\rr}(v_{\rr},B_{1}(0)).
\end{eqnarray*}
This and the minimality of $v_{\rr}$ imply that $\mathcal{H}_{\rr}(\ti{v},B_{1}(0))= \mathcal{H}_{\rr}(v_{\rr},B_{1}(0))$, so using also the strict convexity of $z\mapsto H_{\rr}(\cdot,z)$ we obtain
\begin{eqnarray}\label{5}
\tilde{v}=v_{\rr} \ \ \mbox{a.e. in} \ \ B_{1}(0)\quad \mbox{and}\quad \lim_{j\to \infty}\mathcal{H}_{\rr}(v_{j},B_{1}(0))=\mathcal{H}_{\rr}(v_{\rr},B_{1}(0)).
\end{eqnarray}
\subsubsection*{Step 3: fractional Caccioppoli inequality} The minimality of $v_{j}$ in Dirichlet class $\ti{v}_{j,\rr}+W^{1,\bar{p}}_{0}(B_{1}(0))$ guarantees the validity of the Euler Lagrange equation
\begin{eqnarray}\label{el}
\int_{B_{1}(0)}\langle\partial H_{j}(x,Dv_{h}),D\varphi\rangle \ \dx=0
\end{eqnarray}
for all $\varphi\in W^{1,\bar{p}}_{0}(B_{1}(0))$. We take any vector $h\in \mathbb{R}^{n}$ so that $\snr{h}\le 2^{-10}$, a cut-off function $\eta\in C^{2}_{c}(B_{1}(0))$ so that
$$
\mathds{1}_{B_{3/4}(0)}\le \eta\le \mathds{1}_{B_{5/6}(0)},\qquad \snr{D\eta}^{2}+\snr{D^{2}\eta}\lesssim 1
$$
and test \eqref{el} with $\varphi:=\tau_{-h}(\eta^{2}\tau_{h}v_{j})$. Exploiting the integration by parts formula for finite difference operators, we obtain
\begin{eqnarray*}
0&=&\int_{B_{1}(0)}\langle\tau_{h}\partial H_{j}(x,Dv_{j}),D(\eta^{2}\tau_{h}v_{j})\rangle \ \dx\nonumber \\
&=&\int_{B_{1}(0)}\eta^{2}\langle\tau_{h}\partial H_{j}(x,Dv_{j}),\tau_{h}Dv_{j}\rangle \ \dx\nonumber \\
&+&2\int_{B_{1}(0)}\eta \tau_{h}v_{j}\langle\tau_{h}\partial H_{j}(x,Dv_{j}),D\eta\rangle \ \dx=:\mbox{(I)}+\mbox{(II)}.
\end{eqnarray*}
Let us introduce quantities
\begin{eqnarray*}
\delta:=\min_{\nu\in \texttt{I}_{\kk}}\alpha_{\nu},\qquad A_{\nu,j}:=\left(\nr{a_{\nu,\rr}}_{L^{\infty}(B_{1}(0))}^{2}+[a_{\nu,\rr}]_{0,\alpha;B_{1}(0)}^{2}+(\sigma_{j}^{\nu})^{2}\right)^{\frac{1}{2p_{\nu}-p}}
\end{eqnarray*}
and set for $m\in \left\{p,p_{1},\cdots,p_{\kk}\right\}$
\begin{flalign*}
\mathcal{D}(h):=\left[\snr{Dv_{j}(x+h)}^{2}+\snr{Dv_{j}(x)}^{2}\right],\qquad \mathcal{I}_{m}(h):=\int_{0}^{1}\snr{Dv_{j}(x+\lambda h)}^{m-2}Dv_{j}(x+\lambda h) \ \d\lambda.
\end{flalign*}
Notice that there is no loss of generality in assuming that $\mathcal{D}(h)>0$, otherwise both terms $\mbox{(I)}$-$\mbox{(II)}$ identically vanish. Moreover, consider a nonnegative, radially symmetric mollifier $\phi\in C^{\infty}(B_{1}(0))$, so that $\nr{\phi}_{L^{1}(0)}=1$, let $\phi_{\snr{h}}:=\snr{h}^{-n}\phi(x/\snr{h})$ and regularize for all $\nu\in \texttt{I}_{\kk}$ coefficient $a_{\nu,\rr}(\cdot)$ as done in \cite[Section 5]{comi} via convolution against $\{\phi_{\snr{h}}\}_{\snr{h}>0}$ thus getting $a_{\snr{h}}^{\nu}:=a_{\nu,\rr}*\phi_{\snr{h}}\in C^{\infty}(B_{7/8}(0))$. The newly defined coefficient has the following features:
\begin{eqnarray}\label{ah}
\begin{cases}
\ \nr{a^{\nu}_{\snr{h}}}_{L^{\infty}(B_{7/8}(0))}\le \nr{a_{\nu,\rr}}_{L^{\infty}(B_{1}(0))}\\
\ \snr{a_{\snr{h}}^{\nu}(x)-a_{\nu,\rr}(x)}\le [a_{\nu,\rr}]_{0,\alpha_{\nu};B_{1}(0)}\snr{h}^{\alpha_{\nu}}\\
\ \snr{Da_{\snr{h}}^{\nu}}\le c[a_{\nu,\rr}]_{0,\alpha_{\nu};B_{1}(0)}\snr{h}^{\alpha_{\nu}-1}
\end{cases}\qquad \mbox{for all} \ \ x\in B_{7/8}(0),
\end{eqnarray}
with $c\equiv c(n)$. This will be helpful in a few lines. Finally, we record that whenever $\gamma>1$ and $G\in L^{\frac{\gamma}{\gamma-1}}(B_{1}(0),\mathbb{R}^{n})$, $F\in W^{1,\gamma}_{0}(B_{5/6}(0),\mathbb{R}^{n})$ and $\snr{h}\le 2^{-10}$ it is
\begin{eqnarray}\label{6}
\int_{B_{1}(0)}\langle\tau_{h}G,F\rangle \ \dx=-\snr{h}\int_{B_{1}(0)}\int_{0}^{1}\langle G(x+\lambda h),\partial_{h/\snr{h}}F\rangle \ \d\lambda \ \dx,
\end{eqnarray}
see \cite[(5.29)]{comi}. Now we are ready to estimate terms $\mbox{(I)}$-$\mbox{(II)}$. Notice that
\begin{eqnarray*}
\mbox{(I)}&=&
p\int_{B_{1}(0)}\eta^{2}\langle\tau_{h}(\snr{Dv_{j}}^{p-2}Dv_{j}),\tau_{h}Dv_{j}\rangle \ \dx\nonumber \\
&+&\sum_{\nu=1}^{\kk}p_{\nu}\int_{B_{1}(0)}\eta^{2}(a_{\nu,\rr}(x)+\sigma^{\nu}_{j})\langle\tau_{h}(\snr{Dv_{j}}^{p_{\nu}-2}Dv_{j}),\tau_{h}Dv_{j}\rangle \ \dx\nonumber \\
&+&\sum_{\nu=1}^{\kk}p_{\nu}\int_{B_{1}(0)}\eta^{2}\left(a_{\nu,\rr}(x+h)-a_{\nu,\rr}(x)\right)\langle \snr{Dv_{j}(x+h)}^{p_{\nu}-2}Dv_{j}(x+h),\tau_{h}Dv_{j}\rangle \ \dx\nonumber \\
&=&\mbox{(I)}_{1}+\mbox{(I)}_{2}+\mbox{(I)}_{3}.
\end{eqnarray*}
Via standard monotonicity properties and Lemma \ref{l6} we bound
\begin{eqnarray*}
\mbox{(I)}_{1}+\mbox{(I)}_{2}&\stackrel{\eqref{vpvp}}{\ge}&c\int_{B_{1}(0)}\eta^{2}\snr{\tau_{h}V_{p}(Dv_{j})}^{2} \ \dx\nonumber \\
&+&c\sum_{\nu=1}^{\kk}\int_{B_{1}(0)}\eta^{2}(a_{\nu,\rr}(x)+\sigma^{\nu}_{j})\snr{\tau_{h}V_{p_{\nu}}(Dv_{j})}^{2} \ \dx,
\end{eqnarray*}
with $c\equiv c(n,p,p_{1},\cdots,p_{\kk})$, while by Young inequality and standard properties of translation operators we have
\begin{eqnarray*}
\snr{\mbox{(I)}_{3}}&\le&c\sum_{\nu=1}^{\kk}\snr{h}^{\alpha_{\nu}}[a_{\nu,\rr}]_{0,\alpha_{\nu};B_{1}(0)}\int_{B_{1}(0)}\eta^{2}\mathcal{D}(h)^{\frac{p_{\nu}-1}{2}\pm\frac{p-2}{4}}\snr{\tau_{h}Dv_{j}} \ \dx\nonumber \\
&\le&\varepsilon \int_{B_{1}(0)}\eta^{2}\mathcal{D}(h)^{\frac{p-2}{2}}\snr{\tau_{h}V_{p}(Dv_{j})}^{2} \ \dx\nonumber \\
&+&\frac{c}{\varepsilon}\sum_{\nu=1}^{\kk}\snr{h}^{2\alpha_{\nu}}[a_{\nu,\rr}]_{0,\alpha_{\nu};B_{1}(0)}^{2}\int_{B_{1}(0)}\snr{Dv_{j}}^{2p_{\nu}-p} \ \dx\nonumber \\
&\stackrel{\eqref{vpvp}}{\le}&c\varepsilon\int_{B_{1}(0)}\eta^{2}\snr{\tau_{h}V_{p}(Dv_{j})}^{2} \ \dx+\frac{c\snr{h}^{2\delta}}{\varepsilon}\sum_{\nu=1}^{\kk}\int_{B_{1}(0)}A_{\nu,j}^{2p_{\nu}-p}\snr{Dv_{j}}^{2p_{\nu}-p} \ \dx,
\end{eqnarray*}
for $c\equiv c(n,p,p_{1},\cdots,p_{\kk},\kk)$. Now let us expand term $\mbox{(II)}$:
\begin{eqnarray*}
\mbox{(II)}&=&2p\int_{B_{1}(0)}\eta\tau_{h}v_{j}\langle \tau_{h}\snr{Dv_{j}}^{p-2}Dv_{j},D\eta\rangle \ \dx\nonumber \\
&+&2\sum_{\nu=1}^{\kk}p_{\nu}\sigma^{\nu}_{j}\int_{B_{1}(0)}\eta\tau_{h}v_{j}\langle \tau_{h}(\snr{Dv_{j}}^{p_{\nu}-2}Dv_{j}),D\eta\rangle  \ \dx\nonumber \\
&+&2\sum_{\nu=1}^{\kk}p_{\nu}\int_{B_{1}(0)}\eta a_{\nu,\rr}(x)\tau_{h}v_{j}\langle \tau_{h}(\snr{Dv_{j}}^{p_{\nu}-2}Dv_{j}),D\eta\rangle  \ \dx\nonumber \\
&+&2\sum_{\nu=1}^{\kk}p_{\nu}\int_{B_{1}(0)}\eta\left(a_{\nu,\rr}(x+h)-a_{\nu,\rr}(x)\right)\tau_{h}v_{j}\langle\tau_{h}(\snr{Dv_{j}}^{p_{\nu}-2}Dv_{j}),D\eta\rangle \ \dx\nonumber \\
&=:&\mbox{(II)}_{1}+\mbox{(II)}_{2}+\mbox{(II)}_{3}+\mbox{(II)}_{4}.
\end{eqnarray*}
Set $\mathds{1}_{p}:=1$ if $p\ge 2$ and $\mathds{1}_{p}=0$ when $p\in (1,2)$ and estimate via Lemma \ref{l6}, Young inequality and H\"older inequality with conjugate exponents $\left(\frac{p}{2},\frac{p}{p-2}\right)$ in the superquadratic case and by \eqref{6}, H\"older inequality with conjugate exponents $\left(\frac{p}{2(p-1)},\frac{p}{2-p}\right)$, Jensen inequality and standard properties of translation operators in the subquadratic case:
\begin{eqnarray*}
\snr{\mbox{(II)}_{1}}&\le&\mathds{1}_{p}\snr{\mbox{(II)}_{1}}+(1-\mathds{1}_{p})\snr{\mbox{(II)}_{1}}\nonumber \\
&\le&\varepsilon\mathds{1}_{p}\int_{B_{1}(0)}\eta^{2}\mathcal{D}(h)^{\frac{p-2}{2}}\snr{\tau_{h}Dv_{j}}^{2} \ \dx+\frac{c\mathds{1}_{p}}{\varepsilon}\int_{B_{1}(0)}\snr{D\eta}^{2}\mathcal{D}(h)^{\frac{p-2}{2}}\snr{\tau_{h}v_{j}}^{2} \ \dx\nonumber \\
&+&c\snr{h}(1-\mathds{1}_{p}) \int_{B_{1}(0)}\snr{\mathcal{I}_{p}(h)}\left[\left(\snr{D\eta}^{2}+\snr{D^{2}\eta}\right)\snr{\tau_{h}v_{j}}+\eta\snr{D\eta}\snr{\tau_{h}Dv_{j}}\right] \ \dx\nonumber \\
&\stackrel{\eqref{vpvp}}{\le}&c\varepsilon\int_{B_{1}(0)}\eta^{2}\snr{\tau_{h}V_{p}(Dv_{j})}^{2} \ \dx+\frac{c\mathds{1}_{p}\snr{h}^{2}}{\varepsilon}\int_{B_{1}(0)}\snr{Dv_{j}}^{p} \ \dx\nonumber \\
&+&c\snr{h}(1-\mathds{1}_{p})\left(\int_{B_{5/6}(0)}\snr{\mathcal{I}_{p}(h)}^{\frac{p}{p-1}} \ \dx\right)^{\frac{p-1}{p}}\left(\int_{B_{5/6}(0)}\snr{\tau_{h}v_{j}}^{p} \ \dx\right)^{1/p}\nonumber \\
&+&c\snr{h}^{2}(1-\mathds{1}_{p})\int_{B_{5/6}(0)}\snr{\mathcal{I}_{p}(h)}^{2}\mathcal{D}(h)^{\frac{2-p}{2}} \ \dx\nonumber \\
&\le&c\varepsilon\int_{B_{1}(0)}\eta^{2}\snr{\tau_{h}V_{p}(Dv_{j})}^{2} \ \dx+\frac{c\snr{h}^{2}}{\varepsilon}\int_{B_{1}(0)}\snr{Dv_{j}}^{p} \ \dx\nonumber \\
&+&c\snr{h}^{2}(1-\mathds{1}_{p})\left(\int_{B_{5/6}(0)}\snr{\mathcal{I}_{p}(h)}^{\frac{p}{p-1}} \ \dx\right)^{\frac{2(p-1)}{p}}\left(\int_{B_{5/6}(0)}\mathcal{D}(h)^{p/2} \ \dx\right)^{\frac{2-p}{p}}\nonumber\\
&\le&c\varepsilon\int_{B_{1}(0)}\eta^{2}\snr{\tau_{h}V_{p}(Dv_{j})}^{2} \ \dx+\frac{c\snr{h}^{2}}{\varepsilon}\int_{B_{1}(0)}\snr{Dv_{j}}^{p} \ \dx,
\end{eqnarray*}
for $c\equiv c(n,p)$. Now we abbreviate
\begin{flalign*}
&\mbox{(II)}_{2}^{\nu}:=2p_{\nu}\sigma^{\nu}_{j}\int_{B_{1}(0)}\eta\tau_{h}v_{j}\langle \tau_{h}(\snr{Dv_{j}}^{p_{\nu}-2}Dv_{j}),D\eta\rangle  \ \dx;\nonumber \\
&\mbox{(II)}_{3}^{\nu}:=2p_{\nu}\int_{B_{1}(0)}\eta a_{\nu,\rr}(x)\tau_{h}v_{j}\langle \tau_{h}(\snr{Dv_{j}}^{p_{\nu}-2}Dv_{j}),D\eta\rangle  \ \dx,
\end{flalign*}
and bound by means of \eqref{vpvp}, Lemma \ref{l6}, Young inequality, H\"older inequality with conjugate exponents $\left(\frac{2p_{\nu}-p}{2},\frac{2p_{\nu}-p}{2p_{\nu}-p-2}\right)$ when $p\ge 2$ and $\left(p,\frac{p}{p-1}\right)$, $\left(2p_{\nu}-p,\frac{2p_{\nu}-p}{2p_{\nu}-p-1}\right)$, $\left(\frac{2p_{\nu}-p}{2(p_{\nu}-1)},\frac{2p_{\nu}-p}{2-p}\right)$ if $p\in (1,2)$ and Jensen inequality,
\begin{flalign}\label{d1}
\snr{\mbox{(II)}_{2}^{\nu}}\le&\mathds{1}_{p}\snr{\mbox{(II)}_{2}^{\nu}}+(1-\mathds{1}_{p})\snr{\mbox{(II)}_{2}^{\nu}}\nonumber \\
\le&\varepsilon\mathds{1}_{p}\int_{B_{1}(0)}\eta^{2}\snr{\tau_{h}V_{p}(Dv_{j})}^{2} \ \dx+\frac{c\mathds{1}_{p}(\sigma^{\nu}_{j})^{2}}{\varepsilon}\int_{B_{1}(0)}\snr{D\eta}^{2}\mathcal{D}(h)^{\frac{2p_{\nu}-p-2}{2}}\snr{\tau_{h}v_{j}}^{2} \ \dx\nonumber \\
+&c\snr{h}(1-\mathds{1}_{p})\sigma^{\nu}_{j}\int_{B_{5/6}(0)}\int_{0}^{1}\snr{Dv_{j}(x+\lambda h)}^{p_{\nu}-1}\snr{\tau_{h}v_{j}} \ \d\lambda \  \dx\nonumber \\
+&c\snr{h}(1-\mathds{1}_{p})\sigma^{\nu}_{j}\int_{B_{1}(0)}\eta\snr{D\eta}\snr{\mathcal{I}_{p_{\nu}}(h)}\snr{\tau_{h}Dv_{j}} \ \dx\nonumber\\
\le&\varepsilon\int_{B_{1}(0)}\eta^{2}\snr{\tau_{h}V_{p}(Dv_{j})}^{2} \ \dx+\frac{c\snr{h}^{2}(\sigma^{\nu}_{j})^{2}}{\varepsilon}\int_{B_{1}(0)}\snr{Dv_{j}}^{2p_{\nu}-p} \ \dx\nonumber \\
+&c\snr{h}(1-\mathds{1}_{p})\int_{B_{5/6}(0)}\int_{0}^{1}\left[(\sigma^{\nu}_{j})^{2}\snr{Dv_{j}(x+\lambda h)}^{2p_{\nu}-p-1}+\snr{Dv_{j}(x+\lambda h)}^{p-1}\right]\snr{\tau_{h}v_{j}} \ \dx\nonumber \\
+&c\snr{h}^{2}(1-\mathds{1}_{p})(\sigma^{\nu}_{j})^{2}\left(\int_{B_{5/6}(0)}\snr{\mathcal{I}_{p_{\nu}}(h)}^{\frac{2p_{\nu}-p}{p_{\nu}-1}} \ \dx\right)^{\frac{2(p_{\nu}-1)}{2p_{\nu}-p}}\left(\int_{B_{5/6}(0)}\mathcal{D}(h)^{\frac{2p_{\nu}-p}{2}} \ \dx\right)^{\frac{2-p}{2p_{\nu}-p}}\nonumber \\
\le&\varepsilon\int_{B_{1}(0)}\eta^{2}\snr{\tau_{h}V_{p}(Dv_{j})}^{2} \ \dx+\frac{c\snr{h}^{2}}{\varepsilon}\int_{B_{1}(0)}A_{\nu,j}^{2p_{\nu}-p}\snr{Dv_{j}}^{2p_{\nu}-p} \ \dx ,
\end{flalign}
with $c\equiv c(n,p,p_{\nu})$. Summing the above inequality for $\nu\in \texttt{I}_{\kk}$ we can conclude with
\begin{eqnarray*}
\snr{\mbox{(II)}_{2}}&\le&\sum_{\nu=1}^{\kk}\snr{\mbox{(II)}_{2}^{\nu}}\nonumber \\
&\le&\kk\varepsilon\int_{B_{1}(0)}\eta^{2}\snr{\tau_{h}V_{p}(Dv_{j})}^{2} \ \dx+\frac{c\snr{h}^{2}}{\varepsilon}\sum_{\nu=1}^{\kk}\int_{B_{1}(0)}A_{\nu,j}^{2p_{\nu}-p}\snr{Dv_{j}}^{2p_{\nu}-p} \ \dx,
\end{eqnarray*}
for $c\equiv c(n,p,p_{1},\cdots,p_{\kk},\kk)$. In a similar way, we control
\begin{eqnarray}\label{d2}
\snr{\mbox{(II)}^{\nu}_{3}}&\le&\mathds{1}_{p}\snr{\mbox{(II)}^{\nu}_{3}}+(1-\mathds{1}_{p})\snr{\mbox{(II)}^{\nu}_{3}}\nonumber \\
&\stackrel{\eqref{ah}_{1,2}}{\le}&c\varepsilon\int_{B_{1}(0)}\eta^{2}\snr{\tau_{h}V_{p}(Dv_{j})}^{2} \ \dx\nonumber \\
&+&\frac{c\mathds{1}_{p}\nr{a_{\nu,\rr}}_{L^{\infty}(B_{1}(0))}^{2}}{\varepsilon}\int_{B_{1}(0)}\snr{D\eta}^{2}\mathcal{D}(h)^{p_{\nu}-2-\frac{p-2}{2}}\snr{\tau_{h}v_{j}}^{2} \ \dx\nonumber \\
&+&c(1-\mathds{1}_{p})\snr{h}^{\alpha_{\nu}}\int_{B_{5/6}(0)}\left[[a_{\nu,\rr}]_{0,\alpha_{\nu};B_{1}(0)}^{2}\mathcal{D}(h)^{\frac{2p_{\nu}-p-1}{2}}+\mathcal{D}(h)^{\frac{p-1}{2}}\right]\snr{\tau_{h}v_{j}} \ \dx\nonumber \\
&+&c(1-\mathds{1}_{p})\snr{h}\nr{a_{\nu,\rr}}^{2}_{L^{\infty}(B_{1}(0))}\int_{B_{5/6}(0)}\int_{0}^{1}\snr{Dv_{j}(x+\lambda h)}^{2p_{\nu}-p-1}\snr{\tau_{h}v_{j}} \ \dx\nonumber \\
&+&c(1-\mathds{1}_{p})\snr{h}\int_{B_{5/6}(0)}\int_{0}^{1}\snr{Dv_{j}(x+\lambda h)}^{p-1}\snr{\tau_{h}v_{j}} \ \dx\nonumber \\
&+&c(1-\mathds{1}_{p})\snr{h}^{2}\nr{a_{\nu,\rr}}^{2}_{L^{\infty}(B_{1}(0))}\left(\int_{B_{5/6}(0)}\snr{\mathcal{I}_{p_{\nu}}(j)}^{\frac{2p_{\nu}-p}{p_{\nu}-1}} \ \dx\right)^{\frac{2(p_{\nu}-1)}{2p_{\nu}-p}}\nonumber \\
&\cdot&\left(\int_{B_{5/6}(0)}\mathcal{D}(h)^{\frac{2p_{\nu}-p}{2}} \ \dx\right)^{\frac{2-p}{2p_{\nu}-p}}+c\snr{h}\int_{B_{5/6}(0)}\snr{\mathcal{I}_{p_{\nu}}(h)}\snr{\tau_{h}v_{j}}\snr{Da^{\nu}_{\snr{h}}} \ \dx\nonumber \\
&\stackrel{\eqref{ah}_{3}}{\le}&c\varepsilon\int_{B_{1}(0)}\eta^{2}\snr{\tau_{h}V_{p}(Dv_{j})}^{2} \ \dx+c\snr{h}^{1+\alpha_{\nu}}\int_{B_{1}(0)}\snr{Dv_{j}}^{p} \ \dx\nonumber \\
&+&\frac{c}{\varepsilon}\snr{h}^{1+\alpha_{\nu}}\left(\nr{a_{\nu,\rr}}_{L^{\infty}(B_{1}(0))}^{2}+[a_{\nu,\rr}]_{0,\alpha_{\nu};B_{1}(0)}^{2}\right)\int_{B_{1}(0)}\snr{Dv_{j}}^{2p_{\nu}-p} \ \dx\nonumber \\
&+&c\snr{h}^{\alpha_{\nu}}\int_{B_{5/6}(0)}\int_{0}^{1}[a_{\nu,\rr}]_{0,\alpha_{\nu};B_{1}(0)}^{2}\snr{Dv_{j}(x+\lambda h)}^{2p_{\nu}-1}\snr{\tau_{h}v_{j}} \ \d\lambda \ \dx\nonumber \\
&+&c\snr{h}^{\alpha_{\nu}}\int_{B_{5/6}(0)}\int_{0}^{1}\snr{Dv_{j}(x+\lambda h)}^{p-1}\snr{\tau_{h}v_{j}} \ \d\lambda \ \dx\nonumber \\
&\le&c\varepsilon\int_{B_{1}(0)}\eta^{2}\snr{\tau_{h}V_{p}(Dv_{j})}^{2} \ \dx+c\snr{h}^{1+\alpha_{\nu}}\int_{B_{1}(0)}\snr{Dv_{j}}^{p} \ \dx\nonumber \\
&+&\frac{c\snr{h}^{1+\alpha_{\nu}}}{\varepsilon}\int_{B_{1}(0)}\snr{Dv_{j}}^{2p_{\nu}-p} \ \dx,
\end{eqnarray}
with $c\equiv c(n,p,p_{\nu})$. Summing the inequalities in the previous display we obtain
\begin{eqnarray*}
\snr{\mbox{(II)}_{3}}&\le&\sum_{\nu=1}^{\kk}\snr{\mbox{(II)}^{\nu}_{3}}\nonumber \\
&\le&c\varepsilon\int_{B_{1}(0)}\eta^{2}\snr{\tau_{h}V_{p}(Dv_{j})}^{2} \ \dx+c\snr{h}^{2\delta}\int_{B_{1}(0)}\snr{Dv_{j}}^{p} \ \dx\nonumber \\
&+&\frac{c\snr{h}^{2\delta}}{\varepsilon}\sum_{\nu=1}^{\kk}\int_{B_{1}(0)}A_{\nu,j}^{2p_{\nu}-p}\snr{Dv_{j}}^{2p_{\nu}-p} \ \dx,
\end{eqnarray*}
for $c\equiv c(n,p,p_{1},\cdots,p_{\kk},\kk)$. We stress that when dealing with terms $\mbox{(II)}^{\nu}_{2}$-$\mbox{(II)}^{\nu}_{3}$ we assumed that $\snr{Dv_{j}(x+\lambda h)}>0$ in the sixth and in the fifth and sixth line of displays \eqref{d1}-\eqref{d2} respectively. There is no loss of generality in this as otherwise the integrals in such lines would vanish identically. Concerning term $\mbox{(II)}_{4}$, we have
\begin{eqnarray*}
\snr{\mbox{(II)}_{4}}&\le&c\sum_{\nu=1}^{\kk}[a_{\nu,\rr}]_{0,\alpha_{\nu};B_{1}(0)}\snr{h}^{\alpha_{\nu}}\int_{B_{5/6}(0)}\snr{\tau_{h}v_{j}}\mathcal{D}(h)^{\frac{p_{\nu}-1}{2}} \ \dx\nonumber \\
&\le&c\sum_{\nu=1}^{\kk}\snr{h}^{\alpha_{\nu}}\int_{B_{1}(0)}\left[[a_{\nu,\rr}]_{0,\alpha_{\nu};B_{1}(0)}^{2}\mathcal{D}(h)^{\frac{2p_{\nu}-p-1}{2}}+\mathcal{D}(h)^{\frac{p-1}{2}}\right]\snr{\tau_{h}v_{j}} \ \dx\nonumber \\
&\le&c\snr{h}^{2\delta}\int_{B_{1}(0)}\snr{Dv_{j}}^{p} \ \dx+c\snr{h}^{2\delta}\sum_{\nu=1}^{\kk}\int_{B_{1}(0)}A_{\nu,j}^{2p_{\nu}-p}\snr{Dv_{j}}^{2p_{\nu}-p} \ \dx,
\end{eqnarray*}
with $c\equiv c(n,p,p_{1},\cdots,p_{\kk},\kk)$. Combining the content of all the previous displays and suitably reducing the size of $\varepsilon>0$ we obtain
\begin{eqnarray}\label{7}
\int_{B_{1}(0)}\eta^{2}\snr{\tau_{h}V_{p}(Dv_{j})}^{2} \ \dx &\le&c\snr{h}^{2\delta}\int_{B_{1}(0)}\snr{Dv_{j}}^{p} \ \dx\nonumber\\
&+&c\snr{h}^{2\delta}\sum_{\nu=1}^{\kk}\int_{B_{1}(0)}A_{\nu,j}^{2p_{\nu}-p}\snr{Dv_{j}}^{2p_{\nu}-p} \ \dx,
\end{eqnarray}
for $c\equiv c(n,p,p_{1},\cdots,p_{\kk},\kk)$. At this stage we treat separately the superquadratic case $p\ge 2$ and the subquadratic one $p\in (1,2)$.
\subsubsection*{Step 4: Higher integrability via interpolation - $p\ge 2$} 
From \eqref{vpvp} and \eqref{7} we obtain
\begin{eqnarray*}
\int_{B_{3/4}(0)}\snr{\tau_{h}Dv_{j}}^{p} \ \dx\le c \snr{h}^{2\delta}\int_{B_{1}(0)}\left[\snr{Dv_{j}}^{p}+\sum_{\nu=1}^{\kk}A_{\nu,j}^{2p_{\nu}-p}\snr{Dv_{j}}^{2p_{\nu}-p} \right]\ \dx,
\end{eqnarray*}
with $c\equiv c(n,p,p_{1},\cdots,p_{\kk},\kk)$, so we apply Lemma \ref{l2} to deduce that $Dv_{j}\in W^{s/p,p}(B_{2/3}(0),\mathbb{R}^{n})$ for all $s\in (0,2\delta)$, for simplicity choose $s=\delta$, with
\begin{eqnarray}
\nr{Dv_{j}}_{W^{s/p,p}(B_{2/3}(0))}&\le& c\left[\nr{Dv_{j}}_{L^{p}(B_{1}(0))}+\sum_{\nu=1}^{\kk}\nr{A_{\nu,j}Dv_{j}}_{L^{2p_{\nu}-p}(B_{1}(0))}^{\frac{2p_{\nu}-p}{p}}\right]
\label{8},
\end{eqnarray}
with $c\equiv c(\texttt{data}_{0},\kk)$. Recall that functional $\mathcal{H}_{j}(\cdot)$ is of multi-phase type so Theorem \ref{mor} applies and $v_{j}\in C^{0,\gamma_{0}}(B_{2/3}(0))$ for all $\gamma_{0}\in (0,1)$ therefore for any $0<\chi<\gamma_{0}<1$ and $q\ge 1$ it is 
\begin{eqnarray}\label{9}
[v_{j}]_{\chi,q;B_{2/3}(0)}\le \frac{c[v_{j}]_{0,\gamma_{0};B_{2/3}(0)}}{(q(\gamma_{0}-\chi))^{1/q}}\stackrel{\eqref{hh}}{<}\infty,
\end{eqnarray}
with $c\equiv c(n)$ so we can apply Lemma \ref{fraclem} to get
\begin{eqnarray}
\nr{Dv_{j}}_{L^{t}(B_{1/2}(0))}&\le& c [v_{j}]_{\chi,q;B_{2/3}(0)}^{\theta_{1}}\nr{Dv_{j}}_{W^{s/p,p}(B_{2/3}(0))}^{1-\theta_{1}}\nonumber \\
&\stackrel{\eqref{5}}{\le}&c[v_{j}]_{\chi,q;B_{2/3}(0)}^{\theta_{1}}\left[\nr{Dv_{j}}_{L^{p}(B_{1}(0))}^{1-\theta_{1}}+\sum_{\nu=1}^{\kk}\nr{A_{\nu,j}Dv_{j}}_{L^{2p_{\nu}-p}(B_{1}(0))}^{\frac{(2p_{\nu}-p)(1-\theta_{1})}{p}}\right]
\label{10}
\end{eqnarray}
where $\chi<\gamma_{0}\in (0,1)$, $q>p$ are arbitrary numbers, $c\equiv c(\texttt{data}_{0},\chi,q,\theta_{1},t)$ and it is
\begin{eqnarray}\label{inter}
1=\theta_{1}\chi+(1-\theta_{1})\left(1+s/p\right)\qquad \mbox{and}\qquad \frac{1}{t}=\frac{\theta_{1}}{q}+\frac{1-\theta_{1}}{p},
\end{eqnarray}
which in turn yields that
\begin{eqnarray}\label{11}
\left\{
\begin{array}{c}
\displaystyle 
\ \theta_{1}\equiv \theta_{1}(\chi)=\frac{s/p}{1-\chi+s/p} \ \Longrightarrow \ 1-\theta_{1}=\frac{1-\chi}{1-\chi+s/p}\\ [17 pt]
\displaystyle  \ t\equiv t(q,\chi):=\frac{qp}{p\theta_{1}+q(1-\theta_{1})}=\frac{q(p(1-\chi)+s)}{s+q(1-\chi)},
\end{array}
\right.
\end{eqnarray}
We stress that $\theta_{1}\equiv \theta_{1}(p,\alpha_{1},\cdots,\alpha_{\kk},\chi)$ is increasing with respect to $\chi$ and exponent $t\equiv t(p,q,\alpha,\alpha_{1},\cdots,\alpha_{\kk},\chi)$ is increasing with respect to both, $\chi$ and $q$. Next, we fix $\tau_{1},\tau_{2}\in [1/2,2/3]$, $\tau_{1}<\tau_{2}$ and, following \cite[Section 3.6]{comib} we set $\sigma:=(\tau_{2}-\tau_{1})/4$ and, for a finite $\mathcal{J}\subset \N$, take a covering of $B_{\tau_{1}}(0)$ with a collection of balls $\{B_{\sigma/2}(y_{\iota})\}_{\iota\in \mathcal{J}}$ made by $\snr{\mathcal{J}}=c(n)(\tau_{2}-\tau_{1})^{-n}$ balls so that $y_{\iota}\in B_{\tau_{1}}(0)$ for all $\iota\in \mathcal{J}$. Notice that such a covering can be chosen in such a way that the finite intersection property is satisfied, in the sense that each doubled ball $B_{\sigma}(y_{\iota})$ intersects at most $8^{n}$ of other doubled balls from the same family. We further scale $v_{j}$ on every ball $B_{\sigma}(y_{\iota})$ by defining $v_{\iota}(x):=\sigma^{-1}(v_{j}(y_{\iota}+\sigma x)-(v_{j})_{B_{\sigma}(y_{\iota})})$, $a_{\iota}^{\nu}(x):=a_{\nu,\rr}(y_{\iota}+\sigma x)$ and $H_{\iota}(x,z):=\left[\snr{z}^{p}+\sum_{\nu=1}^{\kk}a_{\iota}^{\nu}(x)\snr{z}^{p_{\nu}}\right]$. Since $v_{j}$ is the solution of \eqref{pda} and therefore it is a local minimizer of functional $\mathcal{H}_{j}(\cdot)$ on $B_{1}(0)$, it is easy to see that $v_{\iota}$ minimizes functional
\begin{eqnarray*}
W^{1,\bar{p}}(B_{1}(0))\ni w\mapsto \min \int_{B_{1}(0)}H_{\iota}(x,Dw) \ \dx,
\end{eqnarray*}
and, keeping \eqref{d3} in mind, we see that \eqref{7} holds for $v_{\iota}$ as well. Recalling that $$[v_{\iota}]_{\chi,q;B_{2/3}(0)}=\sigma^{\chi-1-n/q}[v_{j}]_{\chi,q;B_{2\sigma/3}(y_{\iota})},$$ we can scale back to $v_{j}$ for getting
\begin{eqnarray}\label{12}
\int_{B_{\sigma/2}(y_{\iota})}\snr{Dv_{j}}^{t} \ \dx &\le& c\sigma^{t\theta_{1}(\chi-1)+n\left(1-\frac{t\theta_{1}}{q}-\frac{(1-\theta_{1})t}{p}\right)}[v_{j}]_{\chi,q;B_{2\sigma/3}(y_{\iota})}^{\theta_{1}t}\nonumber \\
&\cdot&\left(\int_{B_{\sigma}(y_{\iota})}\left[\snr{Dv_{\rr}}^{p}+\sum_{\nu=1}^{\kk}A_{\nu,j}^{2p_{\nu}-p}\snr{Dv_{j}}^{2p_{\nu}-p}\right] \ \dx\right)^{\frac{(1-\theta_{1})t}{p}}\nonumber \\
&\stackrel{\eqref{inter}_{2}}{\le}&\frac{c[v_{j}]_{\chi,q;B_{2\sigma/3}(y_{\iota})}^{\theta_{1}t}}{\sigma^{t\theta_{1}(1-\chi)}}\nonumber \\
&\cdot&\left(\int_{B_{\sigma}(y_{\iota})}\left[\snr{Dv_{\rr}}^{p}+\sum_{\nu=1}^{\kk}A_{\nu,j}^{2p_{\nu}-p}\snr{Dv_{j}}^{2p_{\nu}-p}\right] \ \dx\right)^{\frac{(1-\theta_{1})t}{p}},
\end{eqnarray}
where it is $c\equiv c(\texttt{data}_{0},\chi,q,\theta_{1},t)$ and we also used that
\begin{eqnarray*}
\begin{cases}
\ \nr{a_{\iota}^{\nu}}_{L^{\infty}(B_{1}(0))}=\nr{a_{\nu,\rr}}_{L^{\infty}(B_{\sigma}(y_{\iota}))},\\
\ [a_{\iota}^{\nu}]_{0,\alpha_{\nu};B_{1}(0)}=\sigma^{\alpha_{\nu}}[a_{\nu,\rr}]_{0,\alpha_{\nu};B_{\sigma}(y_{\iota})}
\end{cases}\quad \mbox{for all} \ \ \nu\in \texttt{I}_{\kk},
\end{eqnarray*}
which yields that
\begin{eqnarray*}
\nr{a^{\nu}_{\iota}}_{L^{\infty}(B_{1}(0))}^{2}+[a^{\nu}_{\iota}]_{0,\alpha_{\nu};B_{1}(0)}^{2}+(\sigma^{\nu}_{j})^{2}\le A_{\nu,j}^{2p_{\nu}-p}.
\end{eqnarray*}
Summing \eqref{12} for $\iota \in \mathcal{J}$ and using the discrete H\"older inequality $\left(\frac{q}{\theta_{1}t},\frac{p}{t(1-\theta_{1})}\right)$ (legal by means of $\eqref{inter}_{2}$), we obtain
\begin{eqnarray*}
\int_{B_{\tau_{1}}(0)}\snr{Dv_{j}}^{t} \ \dx &\le& \sum_{\iota\in \mathcal{J}}\int_{B_{\sigma/2}(y_{\iota})}\snr{Dv_{j}}^{t} \ \dx\nonumber \\
&\le&\frac{c}{\sigma^{t\theta_{1}(1-\chi)}}\sum_{\iota\in \mathcal{J}}[v_{j}]_{\chi,q;B_{2\sigma/3}(y_{\iota})}^{\theta_{1}t}\nonumber \\
&\cdot&\left(\int_{B_{\sigma}(y_{\iota})}\left[\snr{Dv_{j}}^{p}+\sum_{\nu=1}^{\kk}A_{\nu,j}^{2p_{\nu}-p}\snr{Dv_{j}}^{2p_{\nu}-p}\right] \ \dx\right)^{\frac{(1-\theta_{1})t}{p}}\nonumber \\
&\le&\frac{c}{\sigma^{t\theta_{1}(1-\chi)}}\left(\sum_{\iota\in \mathcal{J}}[v_{j}]_{\chi,q;B_{2\sigma/3}(y_{\iota})}^{q}\right)^{\frac{\theta_{1}t}{q}}\nonumber \\
&\cdot&\left(\sum_{\iota\in \mathcal{J}}\int_{B_{\sigma}(y_{\iota})}\left[\snr{Dv_{j}}^{p}+\sum_{\nu=1}^{\kk}A_{\nu,j}^{2p_{\nu}-p}\snr{Dv_{j}}^{2p_{\nu}-p}\right] \ \dx\right)^{\frac{(1-\theta_{1})t}{p}}\nonumber \\
&\le&\frac{c[v_{j}]_{\chi,q;B_{2/3}(0)}^{\theta_{1}t}}{(\tau_{2}-\tau_{1})^{t\theta_{1}(1-\chi)}}\left(\int_{B_{\tau_{2}}(0)}\left[\snr{Dv_{j}}^{p}+\sum_{\nu=1}^{\kk}A_{\nu,j}^{2p_{\nu}-p}\snr{Dv_{j}}^{2p_{\nu}-p}\right] \ \dx\right)^{\frac{(1-\theta_{1})t}{p}},
\end{eqnarray*}
for $c\equiv c(\texttt{data}_{0},\chi,q,\theta_{1},t)$. Here, we also used that $B_{\sigma}(y_{\iota})\subset B_{\tau_{2}}(0)\subset B_{2/3}(0)$ and that $\mathbb{R}^{n}\ni \omega\mapsto [v_{\rr}]_{\chi,q;\omega}^{q}$ is superadditive as a set function. All in all, using also \eqref{9} and \eqref{hh} we get
\begin{flalign}\label{13}
\nr{Dv_{j}}_{L^{t}(B_{\tau_{1}}(0))}\le&\frac{c[v_{j}]^{\theta_{1}}_{\chi,q;B_{2/3}(0)}}{(\tau_{2}-\tau_{1})^{\theta_{1}(1-\chi)}}\left[\nr{Dv_{j}}_{L^{p}(B_{\tau_{2}}(0))}^{1-\theta_{1}}+\sum_{\nu=1}^{\kk}\nr{A_{\nu,j}Dv_{j}}_{L^{2p_{\nu}-p}(B_{\tau_{2}}(0))}^{\frac{(2p_{\nu}-p)(1-\theta_{1})}{p}}\right]\nonumber \\
\le &\frac{c[v_{j}]^{\theta_{1}}_{0,\gamma_{0};B_{2/3}(0)}}{(\tau_{2}-\tau_{1})^{\theta_{1}(1-\chi)}}\left[\nr{Dv_{j}}_{L^{p}(B_{\tau_{2}}(0))}^{1-\theta_{1}}+\sum_{\nu=1}^{\kk}\nr{A_{\nu,j}Dv_{j}}_{L^{2p_{\nu}-p}(B_{\tau_{2}}(0))}^{\frac{(2p_{\nu}-p)(1-\theta_{1})}{p}}\right]\nonumber \\
\le &\frac{c\mathcal{H}_{j}(v_{j},B_{1}(0))^{\theta_{1}/p}}{(\tau_{2}-\tau_{1})^{\theta_{1}(1-\chi)}}\left[\nr{Dv_{j}}_{L^{p}(B_{\tau_{2}}(0))}^{1-\theta_{1}}+\sum_{\nu=1}^{\kk}\nr{A_{\nu,j}Dv_{j}}_{L^{2p_{\nu}-p}(B_{\tau_{2}}(0))}^{\frac{(2p_{\nu}-p)(1-\theta_{1})}{p}}\right],
\end{flalign}
with $c\equiv c(\texttt{data},\nr{H(\cdot,Dv)}_{L^{1+\delta_{g}}(B_{\rr}(x_{0}))},\delta_{0},\gamma_{0},\chi,q,\theta_{1},t)$. Now fix any $d>\max_{\nu\in \texttt{I}_{\kk}}2p_{\nu}-p$. A straightforward computation yields the chain of implications:
\begin{eqnarray*}
\chi>1-\frac{s}{2d-p} \ \Longrightarrow \ \theta_{1} >1-\frac{p}{2d} \ \Longrightarrow \ \frac{p}{2}-d(1-\theta_{1})>0,
\end{eqnarray*}
which in turn implies that we can choose a suitable lower bound on $q$ so that
\begin{eqnarray*}
q>2d>\frac{dp\theta_{1}}{p-d(1-\theta_{1})} \ \Longrightarrow \ t>d.
\end{eqnarray*}
This means that in \eqref{10} we can use the interpolation inequalities:
\begin{eqnarray*}
\nr{Dv_{j}}_{L^{2p_{\nu}-p}(B_{\tau_{2}}(0))}\le \nr{Dv_{j}}_{L^{t}(B_{\tau_{2}}(0))}^{1-\lambda_{\nu}}\nr{Dv_{j}}_{L^{p}(B_{\tau_{2}}(0))}^{\lambda_{\nu}},
\end{eqnarray*}
where it is
\begin{flalign*}
\frac{1}{2p_{\nu}-p}=\frac{1-\lambda_{\nu}}{t}+\frac{\lambda_{\nu}}{p} \ \Longrightarrow \ \lambda_{\nu}=\frac{p(t+p-2p_{\nu})}{(2p_{\nu}-p)(t-p)}\quad\mbox{and}\quad 1-\lambda_{\nu}=\frac{2t(p_{\nu}-p)}{(2p_{\nu}-p)(t-p)},
\end{flalign*}
for all $\nu\in \texttt{I}_{\kk}$, to have
\begin{eqnarray}\label{14}
\nr{Dv_{j}}_{L^{t}(B_{\tau_{1}}(0))}&\le&\frac{c}{(\tau_{2}-\tau_{1})^{\theta_{1}(1-\chi)}}\mathcal{H}_{j}(v_{j},B_{1}(0))^{1/p}\nonumber \\
&+&\frac{c\mathcal{H}_{j}(v_{j},B_{1}(0))^{\theta_{1}/p}}{(\tau_{2}-\tau_{1})^{\theta_{1}(1-\chi)}}\sum_{\nu=1}^{\kk}A_{\nu,j}^{\frac{(2p_{\nu}-p)(1-\theta_{1})}{p}}\nr{Dv_{j}}_{L^{t}(B_{\tau_{2}}(0))}^{Y_{\nu}/p}\nr{Dv_{j}}_{L^{p}(B_{\tau_{2}}(0))}^{\frac{(2p_{\nu}-p)(1-\theta_{1})\lambda_{\nu}}{p}},
\end{eqnarray}
where 
\begin{eqnarray}\label{yy}
Y_{\nu}:=(2p_{\nu}-p)(1-\theta_{1})(1-\lambda_{\nu}).
\end{eqnarray}
At this stage, we can fix $q=4d$, notice that
\begin{eqnarray}\label{22}
\chi>\chi_{1}:=\max\left\{1-\frac{s}{2d-p},\max_{\nu\in \texttt{I}_{\kk}}\left(1-\frac{s(4d-p)}{8d(p_{\nu}-p)}\right)\right\} \ \Longrightarrow \ Y_{\nu}/p<1
\end{eqnarray}
for all $\nu\in \texttt{I}_{\kk}$ and that
\begin{flalign}\label{23}
&\chi_{2}:=\max\left\{\chi_{1},\max_{\nu\in \texttt{I}_{\kk}}\left(1-\frac{s\mu p(4d-p)}{(p_{\nu}-p)(2n(4d-p)+8\mu pd)}\right)\right\}<\chi\nonumber \\
&\qquad \qquad \qquad \quad\Longrightarrow \ \max_{\nu\in \texttt{I}_{\kk}}\left( \frac{2n(p_{\nu}-p)(1-\theta_{1})}{p(p-Y_{\nu})}\right)<\mu.
\end{flalign}
From \eqref{22} we see that we can apply Young inequality with conjugate exponents $\left(\frac{p}{Y_{\nu}},\frac{p}{p-Y_{\nu}}\right)$ to get
\begin{eqnarray*}
\nr{Dv_{j}}_{L^{t}(B_{\tau_{1}}(0))}&\le&\frac{1}{16}\nr{Dv_{j}}_{L^{t}(B_{\tau_{2}}(0))}+\frac{c}{(\tau_{2}-\tau_{1})^{\theta_{1}(1-\chi)}}\mathcal{H}_{j}(v_{j},B_{1}(0))^{1/p}\nonumber \\
&+&\sum_{\nu=1}^{\kk}\frac{c\mathcal{H}_{j}(v_{j},B_{1}(0))^{\frac{p\theta_{1}+(2p_{\nu}-p)(1-\theta_{1})\lambda_{\nu}}{p(p-Y_{\nu})}}A_{\nu,j}^{\frac{(2p_{\nu}-p)(1-\theta_{1})}{p-Y_{\nu}}}}{(\tau_{2}-\tau_{1})^{\frac{p\theta_{1}(1-\chi)}{p-Y_{\nu}}}},
\end{eqnarray*}
for $c\equiv c(\texttt{data},\nr{H(\cdot,Dv)}_{L^{1+\delta_{g}}(B_{\rr}(x_{0}))},\mu,d)$. Such a dependency can be justified by the fact that all the parameters coming from Lemma \ref{fraclem} ultimately depend only on $(\texttt{data}_{0},\mu,d)$. The content of the previous display legalizes an application of Lemma \ref{iter}, so we obtain
\begin{eqnarray}\label{15}
\nr{Dv_{j}}_{L^{t}(B_{1/2}(0))}&\le&c\mathcal{H}_{j}(v_{j},B_{1}(0))^{1/p}\nonumber \\
&+&c\sum_{\nu=1}^{\kk}\mathcal{H}_{j}(v_{j},B_{1}(0))^{\frac{p\theta_{1}+(2p_{\nu}-p)(1-\theta_{1})\lambda_{\nu}}{p(p-Y_{\nu})}}A_{\nu,j}^{\frac{(2p_{\nu}-p)(1-\theta_{1})}{p-Y_{\nu}}},
\end{eqnarray}
with $c\equiv c(\texttt{data},\nr{H(\cdot,Dv)}_{L^{1+\delta_{g}}(B_{\rr}(x_{0}))},\mu,d)$.
\subsubsection*{Step 5: Higher integrability via interpolation - 1<p<2}
We jump back to \eqref{7} and apply H\"older and Young inequalities with conjugate exponents $\left(\frac{2}{p},\frac{2}{2-p}\right)$ to get
\begin{eqnarray*}
\int_{B_{1}(0)}\eta^{2}\snr{\tau_{h}Dv_{j}}^{p} \ \dx&\stackrel{\eqref{vpvp}}{\le}&c\left(\int_{B_{1}}\eta^{2}\snr{\tau_{h}V_{p}(Dv_{j})}^{2} \ \dx\right)^{p/2}\left(\int_{B_{1}(0)}\eta^{2}\mathcal{D}(h)^{p/2} \ \dx\right)^{\frac{2-p}{2}}\nonumber \\
&\stackrel{\eqref{7}}{\le}&\snr{h}^{\delta p}\int_{B_{1}(0)}\snr{Dv_{j}}^{p} \ \dx\nonumber \\
&+&c\snr{h}^{\delta p}\left(\sum_{\nu=1}^{\kk}\int_{B_{1}(0)}A_{\nu,j}^{2p_{\nu}-p}\snr{Dv_{j}}^{2p_{\nu}-p} \ \dx\right)^{p/2}\left(\int_{B_{1}(0)}\snr{Dv_{j}}^{p} \ \dx\right)^{\frac{2-p}{2}}\nonumber \\
&\le&c\snr{h}^{\delta p}\int_{B_{1}(0)}\snr{Dv_{j}}^{p} \ \dx+c\snr{h}^{\delta p}\sum_{\nu=1}^{\kk}\int_{B_{1}(0)}A_{\nu,j}^{2p_{\nu}-p}\snr{Dv_{j}}^{2p_{\nu}-p} \ \dx,
\end{eqnarray*}
with $c\equiv c(n,p,p_{1},\cdots,p_{\kk},\kk)$, which by Lemma \ref{l2} yields that $Dv_{j}\in W^{s,p}(B_{2/3}(0),\mathbb{R}^{n})$ for all $s\in (0,\delta)$. At this stage, upon choosing $s=\delta/p$, the procedure remains identical to the one described for the superquadratic case, so \eqref{15} holds also when $p\in (1,2)$.
\subsubsection*{Step 6: Conclusions} Notice that $A_{\nu,j}\to A_{\nu}$ as $j\to \infty$, where
\begin{eqnarray}\label{aaaa}
A_{\nu}:=\left(\nr{a_{\nu,\rr}}_{L^{\infty}(B_{1}(0))}^{2}+[a_{\nu,\rr}]^{2}_{0,\alpha_{\nu};B_{1}(0)}\right)^{\frac{1}{2p_{\nu}-p}}.
\end{eqnarray}
Moreover, we can use \eqref{3}, \eqref{5} and weak lower semicontinuity for passing to the limit in \eqref{15} and obtain
\begin{eqnarray*}
\nr{Dv_{\rr}}_{L^{t}(B_{1/2}(0))}&\le&c\mathcal{H}_{\rr}(v_{\rr},B_{1}(0))^{1/p}\nonumber \\
&+&c\sum_{\nu=1}^{\kk}\mathcal{H}_{\rr}(v_{\rr},B_{1}(0))^{\frac{p\theta_{1}+(2p_{\nu}-p)(1-\theta_{1})\lambda_{\nu}}{p(p-Y_{\nu})}}A_{\nu}^{\frac{(2p_{\nu}-p)(1-\theta_{1})}{p-Y_{\nu}}}.
\end{eqnarray*}
Scaling back to $v$, using H\"older inequality on the left-hand side to control the $L^{d}$-average of $v$ (keep in mind that $t>d$) and setting
\begin{eqnarray}\label{gg}
\left\{
\begin{array}{c}
\displaystyle 
\ \Gamma_{1}^{\nu}:=\frac{p\theta_{1}+(2p_{\nu}-p)(1-\theta_{1})\lambda_{\nu}}{p(p-Y_{\nu})}\\ [17 pt]
\displaystyle  \ \Gamma_{2}^{\nu}:=\frac{(2p_{\nu}-p)(1-\theta_{1})}{p-Y_{\nu}}\\ [17pt]
\displaystyle \Gamma_{\nu}:=\frac{2(1-\theta_{1})}{p-Y_{\nu}},
\end{array}
\right.
\end{eqnarray}
we obtain 
\begin{eqnarray}\label{rhrh}
\left(\mint_{B_{\rr/2}(x_{0})}\snr{Dv}^{d} \ \dx\right)^{1/d}&\le& c\left(\mint_{B_{\rr}(x_{0})}H(x,Dv) \ \dx\right)^{1/p}\nonumber \\
&+&c\sum_{\nu=1}^{\kk}A_{\nu}^{\Gamma^{\nu}_{2}}\left(\mint_{B_{\rr}(x_{0})}H(x,Dv) \ \dx\right)^{\Gamma^{\nu}_{1}},
\end{eqnarray}
with $c\equiv c(\texttt{data},\nr{H(\cdot,Dv)}_{L^{1+\delta_{g}}(B_{\rr}(x_{0}))},\mu,d)$. Now notice that the choice of parameters made in \emph{Step 4} and definitions \eqref{yy}-\eqref{gg} yield that
\begin{eqnarray*}
\Gamma_{1}^{\nu}=\frac{2(p_{\nu}-p)(1-\theta_{1})}{p(p-Y_{\nu})}+\frac{1}{p}\stackrel{\eqref{22}}{>}0,
\end{eqnarray*}
therefore with these expansions \eqref{rhrh} becomes
\begin{eqnarray}\label{20}
\left(\mint_{B_{\rr/2}(x_{0})}\snr{Dv}^{d} \ \dx\right)^{1/d}&\le& c\left(\mint_{B_{\rr}(x_{0})}H(x,Dv) \ \dx\right)^{1/p}\nonumber \\
&+&c\sum_{\nu=1}^{\kk}A_{\nu}^{\Gamma_{2}^{\nu}}\left(\mint_{B_{\rr}(x_{0})}H(x,Dv) \ \dx\right)^{\frac{2(p_{\nu}-p)(1-\theta_{1})}{p(p-Y_{\nu})}+\frac{1}{p}},
\end{eqnarray}
with $c\equiv c(\texttt{data},\nr{H(\cdot,Dv)}_{L^{1+\delta_{g}}},\mu,d)$.
\subsubsection*{Step 7: Degenerate phase}
If $\texttt{deg}_{\textnormal{J}}(B_{\rr}(x_{0}))$ is in force, we first set $\mu=1/2$ to remove it from the dependencies of the constants as it will not have a role in this scenario. Furthermore, \eqref{ab} and a quick computation show that 
\begin{eqnarray}\label{20.1}
\nr{a_{\nu}}_{L^{\infty}(B_{\rr}(x_{0}))}\le 4\rr^{\alpha_{\nu}}[a_{\nu}]_{0,\alpha_{\nu};B_{\rr}(x_{0})}+\inf_{x\in B_{\rr}(x_{0})}a_{\nu}(x),
\end{eqnarray}
so \eqref{20.1} and the definition in \eqref{aaaa} yield that
\begin{eqnarray}\label{21}
A_{\nu}^{2p_{\nu}-p}\le 4J^{2}\rr^{2\alpha_{\nu}} [a_{\nu}]_{0,\alpha_{\nu};B_{\rr}(x_{0})}^{2},
\end{eqnarray}
which means that we can rearrange \eqref{20} as
\begin{eqnarray*}
\left(\mint_{B_{\rr/2}(x_{0})}\snr{Dv}^{d} \ \dx\right)^{1/d}&\le&c\left(\mint_{B_{\rr}(x_{0})}H(x,Dv) \ \dx\right)^{1/p}\nonumber \\
&+&c\sum_{\nu=1}^{\kk}J^{\Gamma_{\nu}}\rr^{\Gamma_{\nu}\left(\alpha_{\nu}-\frac{n(p_{\nu}-p)}{p}\right)}\nr{H(\cdot,Dv)}_{L^{1}(B_{\rr}(x_{0}))}^{\frac{\Gamma_{\nu}(p_{\nu}-p)}{p}}\nonumber \\
&\cdot&\left(\mint_{B_{\rr}(x_{0})}H(x,Dv) \ \dx\right)^{1/p}\nonumber\\
&\stackrel{\eqref{bounds}}{\le}&cJ^{\Gamma}\left(\mint_{B_{\rr}(x_{0})}H(x,Dv) \ \dx\right)^{1/p},
\end{eqnarray*}
where $\Gamma:=\max_{\nu\in \texttt{I}_{\kk}}\Gamma_{\nu}$ and $c\equiv c(\texttt{data},\nr{H(\cdot,Dv)}_{L^{1+\delta_{g}}(B_{\rr}(x_{0}))},d)$.
\subsubsection*{Step 8: Nondegenerate/mixed phase}
Assume that either $\texttt{ndeg}_{\textnormal{J}}(B_{\rr}(x_{0}))$ or $\texttt{mix}_{\textnormal{J}}(B_{\rr}(x_{0}))$ is in force. Keeping \eqref{20.1} in mind, this means that either \eqref{21} never holds or that it is verified only for all those indices belonging to $\texttt{d}$. So it is convenient to replace \eqref{21} with $$A_{\nu}^{2p_{\nu}-p}\le 20\left(\nr{a_{\nu}}_{L^{\infty}(B_{\rr}(x_{0}))}^{2}+[a_{\nu}]_{0,\alpha_{\nu};B_{\rr}(x_{0})}^{2}\right),$$ so we can conclude via \eqref{23} that
\begin{eqnarray*}
\left(\mint_{B_{\rr/2}(x_{0})}\snr{Dv}^{d} \ \dx\right)^{1/d}&\le&c\left(\mint_{B_{\rr}(x_{0})}H(x,Dv) \ \dx\right)^{1/p}\nonumber \\
&+&c\rr^{-\mu}\sum_{\nu=1}^{\kk}\nr{H(\cdot,Dv)}_{L^{1}(B_{\rr}(x_{0}))}^{\frac{(p_{\nu}-p)\Gamma_{\nu}}{p}}\left(\mint_{B_{\rr}(x_{0})}H(x,Dv) \ \dx\right)^{1/p}\nonumber \\
&\le&c\rr^{-\mu}\left(\mint_{B_{\rr}(x_{0})}H(x,Dv) \ \dx\right)^{1/p},
\end{eqnarray*}
with $c\equiv c(\texttt{data},\nr{a_{\nu}}_{L^{\infty}(B_{\rr}(x_{0}))},\nr{H(\cdot,Dv)}_{L^{1+\delta_{g}}(B_{\rr}(x_{0}))},\mu,d)$.
\subsubsection*{Step 9: Dependency of constants and their stability under blow up} In \emph{Step 1} we stressed that the functional $\mathcal{H}_{j}(\cdot)$ preserves the multi-phase structure, therefore all the results listed in Section \ref{rmp} apply. In particular, given that we are working on approximating, rescaled problems, we are interested in studying the stability of the constants appearing in Theorem \ref{mor} when it is applied to the sequence $\{v_{j}\}_{j\in \N}$ solutions to \eqref{pda} with respect to scaling and passage to the limit as $j\to\infty$. As already pointed out in \emph{Step 1}, we notice that by Lemma \ref{geh}, the original local minimizer $v$ of functional $\mathcal{H}(\cdot)$ is locally more integrable, in the sense that whenever $B_{\rr}(x_{0})\Subset \Omega$ is any ball with radius $\rr\in (0,1]$, $v\in W^{1,p(1+\delta_{g})}(B_{\rr}(x_{0}))$ for some $\delta_{g}\equiv \delta_{g}(\texttt{data},\nr{H(\cdot,Dv)}_{L^{1}(B_{\rr}(x_{0}))})$. Such information is directly transferred on the blown up map $v_{\rr}$ defined at the very beginning of \emph{Step 1}, which now satisfies $H_{\rr}(\cdot,Dv_{\rr})\in W^{1,1+\delta_{g}}(B_{1}(0))$, where $\delta_{g}\equiv \delta_{g}(\texttt{data},\nr{H(\cdot,Dv)}_{L^{1}(B_{\rr}(x_{0}))})$ is of course the same higher integrability threshold of $v$. As noticed in \eqref{0r}$_{3}$, $H_{\rr}(\cdot,D\ti{v}_{j,\rr})\to H_{\rr}(\cdot,Dv_{\rr})$ strongly in $L^{1+\delta_{g}}(B_{1}(0))$, so if $j\in \N$ is sufficiently large (and up to relabel there is no loss of generality in assuming that $j\ge 1$) it is by \eqref{1},
\begin{eqnarray}\label{16}
\begin{cases}
\ \nr{H_{j}(\cdot,D\ti{v}_{j,\rr})}_{L^{1}(B_{1}(0))}\le\nr{H_{\rr}(\cdot,Dv_{\rr})}_{L^{1}(B_{1}(0))}+1\\
\ \nr{H_{j}(\cdot,D\ti{v}_{j,\rr})}_{L^{1+\delta_{g}}(B_{1}(0))}\le  \nr{H_{\rr}(\cdot,Dv_{\rr})}_{L^{1+\delta_{g}}(B_{1}(0))}+1.
\end{cases}
\end{eqnarray}
 Looking at $v_{j}$, solution to \eqref{pda}, we see that a global higher integrability result applies by means of Lemma \ref{bgeh} with $\delta_{0}\equiv \delta_{g}$, cf. \eqref{0.2} and, by Remark \ref{inc} the dependency of $c$ from $M_{0}$ is nondecreasing and always appears in the form 
\begin{eqnarray}\label{17}
[a_{\nu,\rr}]_{0,\alpha_{\nu};B_{1}(0)}M_{0}^{\frac{p_{\nu}-p}{p}}\qquad \mbox{for all} \ \ \nu\in \texttt{I}_{\kk},
\end{eqnarray} 
where we have also exploited that
\begin{eqnarray}\label{ccc}
[a_{\nu,\rr}+\sigma^{\nu}_{j}]_{0,\alpha_{\nu};B_{1}(0)}\equiv [a_{\nu,\rr}]_{0,\alpha_{\nu};B_{1}(0)}.
\end{eqnarray}
Precisely, by \eqref{16}$_{1}$ it is $M_{0}:=\nr{H_{\rr}(\cdot,Dv_{\rr})}_{L^{1}(B_{1}(0))}+1$, so scaling \eqref{17} back on $B_{\rr}(x_{0})$, we can conclude that
\begin{eqnarray}\label{19}
[a_{\nu,\rr}]_{0,\alpha_{\nu};B_{1}(0)}M_{0}^{\frac{p_{\nu}-p}{p}}&=&\rr^{\alpha_{\nu}-\frac{n(p_{\nu}-p)}{p}}[a_{\nu}]_{0,\alpha_{\nu};B_{\rr}(x_{0})}\left(\nr{H(\cdot,Dv)}_{L^{1}(B_{\rr}(x_{0}))}+1\right)^{\frac{p_{\nu}-p}{p}}\nonumber \\
&\stackrel{\eqref{bounds}}{\le}&[a_{\nu}]_{0,\alpha_{\nu};B_{\rr}(x_{0})}\left(\nr{H(\cdot,Dv)}_{L^{1}(B_{\rr}(x_{0}))}+1\right)^{\frac{p_{\nu}-p}{p}}.
\end{eqnarray}
Recalling that $c$ is nondecreasing in $M_{0}$, we deduce that 
\begin{eqnarray}\label{18}
c(\texttt{data},M_{0})\stackrel{\eqref{17},\eqref{19}}{\le} c(\texttt{data},\nr{H(\cdot,Dv)}_{L^{1}(B_{\rr}(x_{0}))}).
\end{eqnarray}
The same procedure applies for the constant appearing in the local higher integrability result of Lemma \ref{lgeh} with $M=M_{0}$ as by minimality it is
\begin{eqnarray*}
\nr{H_{j}(\cdot,Dv_{j})}_{L^{1}(B_{1}(0))}\le \nr{H_{j}(\cdot,D\ti{v}_{j,\rr})}_{L^{1}(B_{1}(0))}\stackrel{\eqref{16}_{1}}{\le}\nr{H_{\rr}(\cdot,Dv_{\rr})}_{L^{1}(B_{1}(0))}+1,
\end{eqnarray*}
and the dependencies of the constants from Gehring Lemmas have been fixed. We further stress that, looking at the proof of Gehring Lemmas, \cite[Lemmas 4 and 5]{deoh}, \cite[Chapter 6]{giu} and \cite[Theorem 3 and Proposition 1, Chapter 2]{giamodsou}, we can exploit \eqref{18} to make sure that the higher integrability thresholds $\delta_{g}$ and $\sigma_{g}$ depend ultimately on $(\texttt{data},\nr{H(\cdot,Dv)}_{L^{1}(B_{\rr}(x_{0}))})$. From Remark \ref{inc}, we see also that the all the constants appearing in Theorem \ref{mor} are nondecreasing with respect to $M_{g}$, with the (obvious) choice $M_{g}=\nr{H_{\rr}(\cdot,Dv_{\rr})}_{L^{1+\delta_{g}}(B_{1}(0))}+1$. In fact, Lemma \ref{bgeh} renders
\begin{eqnarray*}
\nr{H_{j}(\cdot,Dv_{j})}_{L^{1+\sigma_{g}}(B_{1}(0))}&\le&c\nr{H_{j}(\cdot,D\ti{v}_{j,\rr})}_{L^{1+\delta_{g}}(B_{1}(0))}\stackrel{\eqref{16}_{2}}{\le}c\nr{H_{\rr}(\cdot,Dv_{\rr})}_{L^{1+\delta_{g}}(B_{1}(0))}+c,
\end{eqnarray*}
for $c\equiv c(\texttt{data},\nr{H(\cdot,Dv)}_{L^{1}(B_{\rr}(x_{0}))})$, cf. \eqref{18}. Again, keeping \eqref{ccc} in mind, from \cite{deoh} we have that this dependency is of the form $[a_{\nu,\rr}]_{0,\alpha_{\nu};B_{1}(0)}M_{g}^{\frac{p_{\nu}-p}{p}}$ for all $\nu\in \texttt{I}_{\kk}$, so scaling back we get 
\begin{eqnarray*}
[a_{\nu,\rr}]_{0,\alpha_{\nu};B_{1}(0)}M_{g}^{\frac{p_{\nu}-p}{p}}&=&\rr^{\alpha_{\nu}-\frac{n(p_{\nu}-p)}{p(1+\delta_{g})}}[a_{\nu}]_{0,\alpha_{\nu};B_{\rr}(x_{0})}\left(\nr{H(\cdot,Dv)}_{L^{1+\delta_{g}}(B_{\rr}(x_{0}))}+1\right)^{\frac{p_{\nu}-p}{p}}\nonumber \\
&\stackrel{\eqref{bounds}}{\le}&[a_{\nu}]_{0,\alpha_{\nu};B_{\rr}(x_{0})}\left(\nr{H(\cdot,Dv)}_{L^{1+\delta_{g}}(B_{\rr}(x_{0}))}+1\right)^{\frac{p_{\nu}-p}{p}},
\end{eqnarray*}
so we can conclude that $c(\texttt{data},M_{g})\le c(\texttt{data},\nr{H(\cdot,Dv)}_{L^{1+\delta_{g}}(B_{\rr}(x_{0}))})$. Moreover, looking carefully to the arguments developed in \cite{deoh}, in addition to those described above, another kind of dependency appears that seems to be dangerous for our blow up procedure. In fact, suitably adapting \cite[Corollary 3]{deoh} to our framework, we have constants that are nondecreasing functions of
\begin{eqnarray}\label{dep1}
\begin{cases}
\ [a_{\nu,\rr}]_{0,\alpha_{\nu};B_{1}(0)}\nr{v_{j}}_{L^{\infty}(B_{5/6}(0))}^{p_{\nu}-p}\ \  \mbox{for all} \ \ \nu\in \texttt{I}_{\kk}\quad &\mbox{if} \ \ p(1+\sigma_{g})\le n\\
\ [a_{\nu,\rr}]_{0,\alpha_{\nu};B_{1}(0)}[v_{j}]_{0,\lambda_{g};B_{5/6}(0)}^{p_{\nu}-p}\ \  \mbox{for all} \ \ \nu\in \texttt{I}_{\kk}\quad &\mbox{if} \ \ p(1+\sigma_{g})>n,
\end{cases}
\end{eqnarray}
where $\sigma_{g}\equiv \sigma_{g}(\texttt{data},\nr{H(\cdot,Dv)}_{L^{1}(B_{\rr}(x_{0}))})$ is the higher integrability threshold given by Lemma \ref{bgeh}, $\lambda_{g}:=1-\frac{n}{p(1+\sigma_{g})}$ is the H\"older continuity exponent given by Morrey's embedding theorem and we also used \eqref{ccc}. Now, if $p(1+\sigma_{g})\le n$, we recall from the proof of \cite[Lemma 6]{deoh} that 
\begin{eqnarray*}
\nr{v_{j}}_{L^{\infty}(B_{5/6}(0))}^{p}&\le& c\mint_{B_{1}(0)}H_{j}(x,v_{j}) \ \dx\nonumber \\
&\le& c\mint_{B_{1}(0)}H_{j}(x,Dv_{j}-D\ti{v}_{j,\rr}) \ \dx+c\mint_{B_{1}(0)}H_{j}(x,\ti{v}_{j,\rr}) \ \dx\nonumber \\
&\le&c\mint_{B_{1}(0)}H_{j}(x,D\ti{v}_{j,\rr}) \ \dx\stackrel{\eqref{1},\eqref{16}_{1}}{\le}c\left(\nr{H_{\rr}(\cdot,Dv_{\rr})}_{L^{1}(B_{1}(0))}+1\right)
\end{eqnarray*}
where $c\equiv c(\texttt{data},\nr{H(\cdot,Dv)}_{L^{1}(B_{\rr}(x_{0}))})$ behaves as described in \eqref{17} so no issues about it arise, see also \cite[proof of Theorem 1.1]{comi}. Here, we also exploited the minimality of $v_{j}$, that by construction it is $(\ti{v}_{j,\rr})_{B_{1}(0)}=0$ and Poincar\'e inequality \eqref{sopo}. This means that scaling back to $B_{\rr}(x_{0})$ in $\eqref{dep1}_{1}$ we have
\begin{eqnarray}\label{24}
[a_{\nu,\rr}]_{0,\alpha_{\nu};B_{1}(0)}\nr{v_{j}}_{L^{\infty}(B_{5/6}(0))}^{p_{\nu}-p}&\le&c[a_{\nu,\rr}]_{0,\alpha_{\nu};B_{1}(0)}\left(\nr{H_{\rr}(\cdot,Dv_{\rr})}_{L^{1}(B_{1}(0))}+1\right)^{\frac{p_{\nu}-p}{p}}\nonumber \\
&=&c\rr^{\alpha_{\nu}-\frac{n(p_{\nu}-p)}{p}}[a_{\nu}]_{0,\alpha_{\nu};B_{\rr}(x_{0})}\left(\nr{H(\cdot,Dv)}_{L^{1}(B_{\rr}(x_{0}))}+1\right)^{\frac{p_{\nu}-p}{p}}\nonumber \\
&\stackrel{\eqref{bounds}}{\le}&c\left(\nr{H(\cdot,Dv)}_{L^{1}(B_{\rr}(x_{0}))}+1\right)^{\frac{p_{\nu}-p}{p}},
\end{eqnarray}
for $c\equiv c(\texttt{data},\nr{H(\cdot,Dv)}_{L^{1}(B_{\rr}(x_{0}))})$ (which, as already mentioned, has been treated in \eqref{18}). On the other hand if $p(1+\sigma_{g})>n$, via Morrey embedding theorem, Lemma \ref{bgeh} and Poincar\'e inequality we have
\begin{eqnarray*}
[v_{j}]_{0,\lambda_{g};B_{5/6}(0)}&\le&c\nr{v_{j}}_{W^{1,p(1+\sigma_{g})}(B_{5/6}(0))}\nonumber \\
&\le&c\nr{Dv_{j}}_{L^{p(1+\sigma_{g})}(B_{5/6}(0))}+c\nr{D\ti{v}_{j,\rr}}_{L^{p(1+\sigma_{g})}(B_{5/6}(0))}+c\nr{\ti{v}_{j,\rr}}_{L^{p(1+\sigma_{g})}(B_{5/6}(0))}\nonumber \\
&\le&c\nr{H_{j}(\cdot,D\ti{v}_{j,\rr})}_{L^{1+\sigma_{g}}(B_{1}(0))}^{1/p}\nonumber \\
&\stackrel{\eqref{16}_{2}}{\le}&c\left(\nr{H_{\rr}(\cdot,Dv_{\rr})}_{L^{1+\delta_{g}}(B_{1}(0))}+1\right)^{1/p}
\end{eqnarray*}
for $c\equiv c(\texttt{data},\nr{H(\cdot,Dv)}_{L^{1}(B_{\rr}(x_{0}))})$ and we also used that $(\ti{v}_{j,\rr})_{B_{1}(0)}=0$. With this last inequality at hand, we can jump back to $\eqref{dep1}_{2}$ and conclude as in \eqref{24}.

\begin{remark}\label{rrr}
\emph{We stress that the constants appearing in \eqref{fh}-\eqref{hf} are nondecreasing with respect to $\nr{H(\cdot,Dv)}_{L^{1+\delta_{g}}(B_{\rr}(x_{0}))}$ and to $J$.}
\end{remark}
\section{Applications to Calder\'on Zygmund estimates}\label{czss}
In this section we provide Calder\'on-Zygmund type estimates for local minimizers of the nonhomogeneous functional $\mathcal{G}(\cdot)$. With \eqref{assg}-\eqref{assf} in force, the definition of minima in this case is the same given in Definition \ref{min} - just replace $\mathcal{H}(\cdot)$ with $\mathcal{G}(\cdot)$ there.
\subsection{Proof of Theorem \ref{czt}}
The outline of the proof of Theorem \ref{czt} is analogous to the one of \cite{bbo,comicz,demicz}, therefore we shall follow the same steps indicated there and point out only the relevant changes.
\subsubsection*{Step 1 - Existence and uniform higher integrability}
Existence and uniqueness for minima of functional $\mathcal{G}(\cdot)$ follows by direct methods under the minimal assumptions $0\le a_{\nu}(\cdot)\in L^{\infty}(\Omega)$ for all $\nu\in \texttt{I}_{\kk}$ and $H(\cdot,F)\in L^{1}(\Omega)$, that are in any case guaranteed by \eqref{ab}, \eqref{assg} and \eqref{assf}, cf. \cite[Remark 1.2]{comicz}. Moreover, a straightforward manipulation of \cite[Theorem 4]{demicz} assures that there is a positive higher integrability threshold $\delta_{\gamma}\equiv \delta_{\gamma}(\texttt{data},\Lambda,\nr{H(\cdot,Du)}_{L^{1}(\ti{\Omega}_{0})})<\gamma-1$ so that
\begin{eqnarray}\label{cz0}
H(\cdot,Du)\in L^{1+\delta_{\gamma}}_{\loc}(\Omega)
\end{eqnarray}
and whenever $B_{\rr}(x_{0})\Subset \Omega$ is a ball with radius $\rr\in (0,1]$ it is
\begin{eqnarray}\label{cz1}
\left(\mint_{B_{\rr/2}(x_{0})}H(x,Du)^{1+\delta} \ \dx \right)^{\frac{1}{1+\delta}}&\le&c\mint_{B_{\rr}(x_{0})}H(x,Du) \ \dx\nonumber \\
&+&c\left(\mint_{B_{\rr}(x_{0})}H(x,F)^{1+\delta} \ \dx\right)^{\frac{1}{1+\delta}}
\end{eqnarray}
for all $\delta\in (0,\delta_{\gamma}]$ with $c\equiv c(\texttt{data},\Lambda,\nr{H(\cdot,Du)}_{L^{1}(B_{\rr}(x_{0}))},\gamma)$.
\subsubsection*{Step 2 - Exit time and covering of level sets}
Let $\Omega_{0}\Subset \ti{\Omega}_{0}\Subset \Omega$ be three open set as in the statement of Theorem \ref{czt} and $B_{r}\Subset \Omega_{0}$ be a ball with radius $r\le r_{*}$, a threshold that will be fixed in a few lines. We recall that \eqref{cz0}-\eqref{cz1} and a standard covering argument render
\begin{flalign}\label{cz2}
\nr{H(\cdot,Du)}_{L^{1+\delta_{\gamma}}(\Omega_{0})}\le c(\texttt{data},\Lambda,\nr{H(\cdot,Du)}_{L^{1}(\ti{\Omega}_{0})},\nr{H(\cdot,F)}_{L^{\gamma}(\ti{\Omega}_{0})},\gamma,\dist(\ti{\Omega}_{0},\partial \Omega)).
\end{flalign}
We apply the exit time and covering argument as in \cite[Theorem 1.1]{comicz}, which in particular yields the collection of balls $\{B_{\iota}\}\equiv \{B_{\rr_{\iota}}(x_{\iota})\}\equiv \{5\ti{B}_{\iota}\}$ as denoted in \cite[(4.9)-(4.11)]{comicz}. All such balls are contained in $B_{r}\Subset \Omega_{0}$. 
\subsubsection*{Step 3 - Comparison, first time}
We construct a first comparison problem. Precisely, we let $v_{\iota}\in u+W^{1,p}_{0}(4B_{\iota})$ be the solution of Dirichlet problem
\begin{eqnarray}\label{pdi}
u+W^{1,p}_{0}(4B_{\iota})\ni w\mapsto \min \mathcal{H}(w,4B_{\iota}),
\end{eqnarray}
whose existence and uniqueness is guaranteed by standard direct methods. By minimality, $v_{\iota}$ satisfied the integral identity
\begin{eqnarray}\label{eli}
0=\mint_{4B_{\iota}}\langle \partial H(x,Dv_{\iota}),D\varphi\rangle \ \dx,
\end{eqnarray}
for all $\varphi\in W^{1,p}_{0}(4B_{\iota})$ so that $H(\cdot,D\varphi)\in L^{1}(4B_{\iota})$. Moreover, by the minimality of $v_{\iota}$ in Dirichlet class $u+W^{1,p}_{0}(4B_{\iota})$, \eqref{cz2}, Lemma \ref{bgeh} and Remark \ref{inc} we have 
\begin{eqnarray}\label{cz3}
\left\{
\begin{array}{c}
\displaystyle 
\ \mint_{4B_{\iota}}H(x,Dv_{\iota}) \ \dx \le \mint_{4B_{\iota}}H(x,Du) \ \dx\\ [17 pt]
\displaystyle  \ \mint_{4B_{\iota}}H(x,Dv_{\iota})^{1+\sigma_{g}} \ \dx \le c\mint_{4B_{\iota}}H(x,Du)^{1+\sigma_{g}} \ \dx,
\end{array}
\right.
\end{eqnarray}
for $c,\sigma_{g}\equiv c,\sigma_{g}(\texttt{data}_{\textnormal{cz}})$ and $\sigma_{g}\in (0,\delta_{\gamma})$. To get this dependency, motivated by \eqref{cz2} and $\eqref{cz3}_{1}$, we choose in Lemma \ref{bgeh} $M_{0}=\nr{H(\cdot,Du)}_{L^{1}(\ti{\Omega}_{0})}$. Moreover, by Theorem \ref{ht} we have that $v_{\iota}\in C^{1,\beta_{0}}_{\loc}(4B_{\iota})$ for some $\beta_{0}\equiv \beta_{0}(\texttt{data}_{0})$ and, according to Theorem \ref{revh}, reverse H\"older inequalities \eqref{fh}-\eqref{hf} hold for all $d\in [1,\infty)$ and any $\mu\in (0,1]$. Extending $u-v_{\iota}\equiv 0$ in $\Omega\setminus 4B_{\iota}$, we see that we can proceed as in \cite[(4.17)]{comicz} to get
\begin{eqnarray}\label{cz4}
\mathcal{V}(Du,Dv_{\iota}):=\mint_{4B_{\iota}}\snr{V_{p}(Du)-V_{p}(Dv_{\iota})}^{2} \ \dx&+&\sum_{\nu=1}^{\kk}\int_{4B_{\iota}}a_{\nu}(x)\snr{V_{p_{\nu}}(Du)-V_{p_{\nu}}(Dv_{\iota})}^{2} \ \dx \nonumber \\
&\le&c\varepsilon\mint_{4B_{\iota}}H(x,Du) \ \dx +\frac{c}{\varepsilon^{\bar{p}}}\mint_{4B_{\iota}}H(x,F) \ \dx,
\end{eqnarray}
for $c\equiv c(n,\Lambda,p,p_{1},\cdots,p_{\kk},\kk)$ and $\bar{p}:=\max_{\nu\in \texttt{I}_{\kk}}p_{\nu}$.
\subsubsection*{Step 4 - Comparison, second time} We define
\begin{eqnarray}\label{++}
a^{+}_{\iota,\nu}:=\sup_{x\in \overline{2B}_{\iota}}a_{\nu}(x)\quad \mbox{for all} \ \ \nu\in \texttt{I}_{\kk}
\end{eqnarray}
and notice that Theorem \ref{ht} yields that $v_{\iota}\in W^{1,\infty}(2B_{\iota})$, therefore setting 
$$
\mathbb{R}^{n}\ni z\mapsto H_{+}(z):=\snr{z}^{p}+\sum_{\nu=1}^{\kk}a^{+}_{\iota,\nu}\snr{z}^{p_{\nu}},
$$
it trivially holds that $H_{+}(Dv_{\iota})\in L^{1}(2B_{\iota})$. This means that we can consider the solution $w_{\iota}\in v_{\iota}+W^{1,p}(2B_{\iota})$ of the second Dirichlet problem
\begin{eqnarray}\label{pd2}
v_{\iota}+W^{1,p}_{0}(2B_{\iota})\ni w\mapsto \min \int_{2B_{\iota}}H_{+}(Dw) \ \dx.
\end{eqnarray}
By minimality, $w_{\iota}$ satisfies
\begin{flalign}\label{el2}
\left\{
\begin{array}{c}
\displaystyle 
\ \mint_{2B_{\iota}}\langle \partial H_{+}(Dw_{\iota}),D\varphi\rangle \ \dx=0\\ [17 pt]
\displaystyle  \ \mint_{2B_{\iota}}H_{+}(Dw_{\iota}) \ \dx\le \mint_{2B_{\iota}}H_{+}(Dv_{\iota}) \ \dx,
\end{array}
\right.
\end{flalign}
and $\eqref{el2}_{1}$ holds for all $\varphi\in W^{1,p}_{0}(2B_{\iota})$ so that $H_{+}(D\varphi)\in L^{1}(2B_{\iota})$. After extending $v_{\iota}-w_{\iota}\equiv 0$ in $\Omega_{0}\setminus 2B_{\iota}$, we see that the function $v_{\iota}-w_{\iota}$ is admissible in both \eqref{eli}-\eqref{el2}; standard monotonicity arguments then yield

\begin{eqnarray}\label{cz5}
\mathcal{V}_{0}(Dv_{\iota},Dw_{\iota})&:=&\mint_{2B_{\iota}}\snr{V_{p}(Dv_{\iota})-V_{p}(Dw_{\iota})}^{2} \ \dx\nonumber\\ &+&\sum_{\nu=1}^{\kk}\mint_{2B_{\iota}}a^{+}_{\iota,\nu}\snr{V_{p_{\nu}}(Dv_{\iota})-V_{p_{\nu}}(Dw_{\iota})}^{2} \ \dx\nonumber \\
&\le&c\mint_{2B_{\iota}}\langle \partial H_{+}(Dv_{\iota})-\partial H_{+}(Dw_{\iota}),Dv_{\iota}-Dw_{\iota}\rangle \ \dx\nonumber \\
&\stackrel{\eqref{eli},\eqref{el2}_{1}}{=}&c\mint_{2B_{\iota}}\langle \partial H_{+}(Dv_{\iota})-\partial H(x,Dv_{\iota}),Dv_{\iota}-Dw_{\iota}\rangle \ \dx\nonumber \\
&\le&c\sum_{\nu=1}^{\kk}\mint_{2B_{\iota}}\snr{a^{+}_{\iota,\nu}-a_{\nu}(x)}\snr{Dv_{\iota}}^{p_{\nu}-1}\snr{Dv_{\iota}-Dw_{\iota}} \ \dx\nonumber \\
&\le&c\sum_{\nu=1}^{\kk}\left(\osc_{ 2B_{\iota}}a_{\nu}\right)\mint_{2B_{\iota}}\snr{Dv_{\iota}}^{p_{\nu}-1}\snr{Dv_{\iota}-Dw_{\iota}} \ \dx\nonumber \\
&=:&\sum_{\nu=1}^{\kk}\mbox{(I)}_{\nu},
\end{eqnarray}
for $c\equiv c(n,p,p_{1},\cdots,p_{\kk},\kk)$. In the following we shall introduce three new positive constants, which may vary from line to line, but will always have the same dependencies:
\begin{itemize}
\item $c_{\texttt{nd}}\equiv c_{\texttt{nd}}(n,\Lambda,p,p_{1},\cdots,p_{\kk},\kk)$;
\item $c_{\texttt{m}}\equiv c_{\texttt{m}}(\texttt{data}_{\textnormal{cz}})$;
\item $c_{\texttt{d}}\equiv c_{\texttt{d}}(\texttt{data},\Lambda,\nr{H(\cdot,Du)}_{L^{1}(\ti{\Omega}_{0})},\nr{H(\cdot,F)}_{L^{\gamma}(\ti{\Omega}_{0})},\gamma,\dist(\ti{\Omega}_{0},\partial \Omega))$.
\end{itemize}
\subsubsection*{Step 5 - Estimates in the nondegenerate phase} Assume that $\texttt{ndeg}_{\textnormal{J}}(2B_{\iota})$ is in force for some $J\ge 4$ that will eventually be fixed as a function of $(n,\Lambda,p,p_{1},\cdots,p_{\kk},\kk)$. In this setting, it is
\begin{eqnarray}\label{cz6}
\osc_{2B_{\iota}}a_{\nu}\le 4\rr_{\iota}^{\alpha_{\nu}}[a_{\nu}]_{0,\alpha_{\nu};2B_{\iota}}\le \frac{4a_{\nu}(x)}{J}\qquad \mbox{for all} \ \ \nu\in \texttt{I}_{\kk}.
\end{eqnarray}
Notice that the very definition of $H_{+}(\cdot)$ and the minimality of $w_{\iota}$ in class $v_{\iota}+W^{1,p}_{0}(2B_{\iota})$ and of $v_{\iota}$ in class $u+W^{1,p}_{0}(4B_{\iota})$ yield that
\begin{eqnarray}\label{cz8}
\mint_{2B_{\iota}}H(x,Dw_{\iota}) \ \dx &\le&\mint_{2B_{\iota}}H_{+}(Dw_{\iota}) \ \dx\nonumber \\
&\le&\mint_{2B_{\iota}}H_{+}(Dv_{\iota}) \ \dx\nonumber \\
&\le&\mint_{2B_{\iota}}H(x,Dv_{\iota}) \ \dx+\sum_{\nu=1}^{\kk}\left(\osc_{2B_{\iota}}a_{\nu}\right)\mint_{2B_{\iota}}\snr{Dv_{\iota}}^{p_{\nu}} \ \dx\nonumber \\
&\stackrel{\eqref{cz6}}{\le}&10\mint_{2B_{\iota}}H(x,Dv_{\iota}) \ \dx \le 4^{2n}\mint_{4B_{\iota}}H(x,Du) \ \dx,
\end{eqnarray}
so we may estimate via H\"older inequality with conjugate exponents $\left(p_{\nu},\frac{p_{\nu}}{p_{\nu}-1}\right)$,
\begin{eqnarray*}
\mbox{(I)}_{\nu}&\stackrel{\eqref{cz6}}{\le}& \frac{c_{\texttt{nd}}}{J}\mint_{2B_{\iota}}a_{\nu}(x)\snr{Dv_{\iota}}^{p_{\nu}-1}\snr{Dv_{\iota}-Dw_{\iota}} \ \dx\nonumber \\
&\le&\frac{c_{\texttt{nd}}}{J}\left(\mint_{2B_{\iota}}a_{\nu}(x)\snr{Dv_{\iota}}^{p_{\nu}} \ \dx\right)^{\frac{p_{\nu}-1}{p_{\nu}}}\left(\mint_{2B_{\iota}}a_{\nu}(x)\left[\snr{Dv_{\iota}}^{p_{\nu}}+\snr{Dw_{\iota}}^{p_{\nu}}\right] \ \dx\right)^{1/p_{\nu}}\nonumber\\
&\stackrel{\eqref{cz8}}{\le}&\frac{c_{\texttt{nd}}}{J}\mint_{4B_{\iota}}H(x,Du) \ \dx,
\end{eqnarray*}
for $c\equiv c(n,p,p_{\nu})$. Summing the content of the above display over $\nu \in \texttt{I}_{\kk}$ we obtain
\begin{eqnarray}\label{cz7}
\sum_{\nu=1}^{\kk}\mbox{(I)}_{\nu}&\le&\frac{c_{\texttt{nd}}}{J}\mint_{4B_{\iota}}H(x,Du) \ \dx.
\end{eqnarray}
\subsubsection*{Step 8 - Estimates in the mixed phase} Now we assume that $\texttt{mix}_{\textnormal{J}}(2B_{\iota})$ holds with $J\ge 4$ still to be fixed, pick any 
\begin{flalign}
\mu\in \left(0,\min_{\nu\in \texttt{I}_{\kk}}\frac{1}{p_{\nu}}\left(\alpha_{\nu}-\frac{n(p_{\nu}-p)}{p(1+\delta_{\gamma})}\right)\right)\stackrel{\eqref{bounds}}{\not =}\emptyset \ \Longrightarrow \ \mu\equiv \mu(\texttt{data},\nr{H(\cdot,Du)}_{L^{1}(\Omega_{0})}),\label{mu}
\end{flalign}
where $\delta_{\gamma}$ is the higher integrability exponent determined in \emph{Step 1} and set $$\sigma_{0}:=\min_{\nu\in \texttt{I}_{\kk}}\left(\alpha_{\nu}-\mu p_{\nu}-\frac{n(p_{\nu}-p)}{p(1+\delta_{\gamma})}\right)\stackrel{\eqref{bounds},\eqref{mu}}{>}0.$$
Keeping in mind that
\begin{eqnarray}\label{osc}
\osc_{2B_{\iota}}a_{\nu}\le 2a^{+}_{\iota,\nu},
\end{eqnarray}
we can proceed as in \cite[Section 6]{bbo} and apply \eqref{hf} with $d=p_{\nu}$ and $\mu$ as in \eqref{mu} to control
\begin{eqnarray*}
\mbox{(I)}_{\nu}&\le&c\rr_{\iota}^{\alpha_{\nu}}\mint_{2B_{\iota}}\snr{Dv_{\iota}}^{p_{\nu}} \ \dx\nonumber \\
&+&\left(\osc_{2B_{\iota}}a_{\nu}\right)\snr{2B_{\iota}}^{-1}\int_{2B_{\iota}\cap\{\snr{Dw_{\iota}}\ge J\snr{Dv_{\iota}}\}}\snr{Dv_{\iota}}^{p_{\nu}-1}\snr{Dw_{\iota}} \ \dx\nonumber \\
&+&\left(\osc_{2B_{\iota}}a_{\nu}\right)\snr{2B_{\iota}}^{-1}\int_{2B_{\iota}\cap\{\snr{Dw_{\iota}}< J\snr{Dv_{\iota}}\}}\snr{Dv_{\iota}}^{p_{\nu}-1}\snr{Dw_{\iota}} \ \dx\nonumber \\
&\stackrel{\eqref{osc}}{\le}&c(1+J)\rr_{\iota}^{\alpha_{\nu}}\mint_{2B_{\iota}}\snr{Dv_{\iota}}^{p_{\nu}} \ \dx+\frac{c}{J^{p-1}}\mint_{2B_{\iota}}a^{+}_{\iota,\nu}\snr{Dw_{\iota}}^{p_{\nu}}\ \dx\nonumber \\
&\stackrel{\eqref{el2}_{2}}{\le}&c(1+J)\rr_{\iota}^{\alpha_{\nu}}\mint_{2B_{\iota}}\snr{Dv_{\iota}}^{p_{\nu}} \ \dx+\frac{c}{J^{p-1}}\mint_{2B_{\iota}}H_{\iota}(Dv_{\iota}) \ \dx\nonumber\\
&\stackrel{\eqref{cz3}_{1}}{\le}&cJ\rr_{\iota}^{\alpha_{\nu}}\mint_{2B_{\iota}}\snr{Dv_{\iota}}^{p_{\nu}} \ \dx+\frac{c}{J^{p-1}}\sum_{m=1}^{\kk}\rr_{\iota}^{\alpha_{m}}\mint_{2B_{\iota}}\snr{Dv_{\iota}}^{p_{m}} \ \dx+\frac{c_{\texttt{nd}}}{J^{p-1}}\mint_{4B_{\iota}}H(x,Du) \ \dx\nonumber\\
&\stackrel{\eqref{hf}}{\le}&c_{\texttt{m}}J\rr_{\iota}^{\alpha_{\nu}-p_{\nu}\mu}\left(\mint_{4B_{\iota}}H(x,Dv_{\iota})^{1+\sigma_{g}} \ \dx\right)^{\frac{p_{\nu}-p}{p(1+\sigma_{g})}}\mint_{4B_{\iota}}H(x,Du) \ \dx\nonumber \\
&+&\frac{c}{J^{p-1}}\sum_{m=1}^{\kk}\rr_{\iota}^{\alpha_{m}-\mu p_{m}}\left(\mint_{4B_{\iota}}H(x,Dv_{\iota})^{1+\sigma_{g}} \ \dx\right)^{\frac{p_{m}-p}{p(1+\sigma_{g})}}\mint_{4B_{\iota}}H(x,Du) \ \dx\nonumber\\
&+&\frac{c_{\texttt{nd}}}{J^{p-1}}\mint_{4B_{\iota}}H(x,Du) \ \dx\nonumber \\
&\stackrel{\eqref{cz3}}{\le}&c_{\texttt{m}}J\rr_{\iota}^{\sigma_{0}}\nr{H(\cdot,Du)}_{L^{1+\delta_{\gamma}}(4B_{\iota})}^{\frac{p_{\nu}-p}{p}}\mint_{4B_{\iota}}H(x,Du) \ \dx+\frac{c_{\texttt{nd}}}{J^{p-1}}\mint_{4B_{\iota}}H(x,Du) \ \dx\nonumber \\
&+&\frac{c_{\texttt{m}}\rr_{\iota}^{\sigma_{0}}}{J^{p-1}}\left(\mint_{4B_{\iota}}H(x,Du) \ \dx\right)\sum_{m=1}^{\kk}\nr{H(\cdot,Du)}_{L^{1+\delta_{\gamma}}(4B_{\iota})}^{\frac{p_{m}-p}{p}}\nonumber\\
&\le&\left(c_{\texttt{m}}J\rr_{\iota}^{\sigma_{0}}+\frac{c_{\texttt{nd}}}{J^{p-1}}\right)\mint_{4B_{\iota}}H(x,Du) \ \dx.
\end{eqnarray*}
We stress that here we also used Remark \ref{rrr} and \eqref{cz3} to determine such a dependency for $c$. Summing the above inequalities over $\nu\in \texttt{I}_{\kk}$ we get
\begin{eqnarray}\label{cz10}
\sum_{\nu=1}^{\kk}\mbox{(I)}_{\nu}&\le&\left(c_{\texttt{m}}J\rr_{\iota}^{\sigma_{0}}+\frac{c_{\texttt{nd}}}{J^{p-1}}\right)\mint_{4B_{\iota}}H(x,Du) \ \dx.
\end{eqnarray}

\subsubsection*{Step 9 - Estimates in the degenerate phase} Finally, we look at the case $\texttt{deg}(2B_{\iota})$. We set
$$
\tau_{0}:=\min_{\nu\in \texttt{I}_{\kk}}\left(\alpha_{\nu}-\frac{n(p_{\nu}-p)}{p(1+\delta_{\gamma})}\right)\stackrel{\eqref{bounds}}{>}0 \ \Longrightarrow \ \tau_{0}\equiv \tau_{0}(\texttt{data},\nr{H(\cdot,Du)}_{L^{1}(\Omega_{0})})
$$
and as done in \emph{Step 8} we estimate
\begin{eqnarray*}
\mbox{(I)}_{\nu}&\le&c\rr_{\iota}^{\alpha_{\nu}}\mint_{2B_{\iota}}\snr{Dv_{\iota}}^{p_{\nu}} \ \dx\nonumber \\
&+&c\snr{2B_{\iota}}^{-1}\left(\osc_{2B_{\iota}}a_{\nu}\right)\int_{2B_{\iota}\cap\{\snr{Dw_{\iota}}\ge J\snr{Dv_{\iota}}\}}\snr{Dv_{\iota}}^{p_{\nu}-1}\snr{Dw_{\iota}} \ \dx\nonumber \\
&+&c\snr{2B_{\iota}}^{-1}\left(\osc_{2B_{\iota}}a_{\nu}\right)\int_{2B_{\iota}\cap\{\snr{Dw_{\iota}}< J\snr{Dv_{\iota}}\}}\snr{Dv_{\iota}}^{p_{\nu}-1}\snr{Dw_{\iota}} \ \dx\nonumber \\
&\stackrel{\eqref{osc}}{\le}&c(1+J)\rr_{\iota}^{\alpha_{\nu}}\mint_{2B_{\iota}}\snr{Dv_{\iota}}^{p_{\nu}} \ \dx+\frac{c}{J^{p-1}}\mint_{2B_{\iota}}a_{\iota,\nu}^{+}\snr{Dw_{\iota}}^{p_{\nu}} \ \dx\nonumber \\
\nonumber \\
&\stackrel{\eqref{el}_{2}}{\le}&cJ\rr_{\iota}^{\alpha_{\nu}}\mint_{2B_{\iota}}\snr{Dv_{\iota}}^{p_{\nu}} \ \dx+\frac{c}{J^{p-1}}\sum_{m=1}^{\kk}\rr_{\iota}^{\alpha_{m}}\mint_{2B_{\iota}}\snr{Dv_{\iota}}^{p_{m}} \ \dx+\frac{c}{J^{p-1}}\mint_{2B_{\iota}}H(x,Dv_{\iota}) \ \dx\nonumber\\
&\stackrel{\eqref{fh},\eqref{cz3}_{1}}{\le}&c_{\texttt{d}}J^{\bar{p}\Gamma+1}\rr^{\alpha_{\nu}}_{\iota}\left(\mint_{4B_{\iota}}H(x,Dv_{\iota})^{1+\sigma_{g}} \ \dx\right)^{\frac{p_{\nu}-p}{p(1+\sigma_{g})}}\mint_{4B_{\iota}}H(x,Dv_{\iota}) \ \dx\nonumber \\
&+&c_{\texttt{d}}J^{\bar{p}\Gamma-p+1}\sum_{m=1}^{\kk}\rr_{\iota}^{\alpha_{m}}\left(\mint_{4B_{\iota}}H(x,Dv_{\iota})^{1+\sigma_{g}} \ \dx\right)^{\frac{p_{m}-p}{p(1+\sigma_{g})}}\mint_{4B_{\iota}}H(x,Dv_{\iota}) \ \dx\nonumber\\
&+&\frac{c_{\texttt{nd}}}{J^{p-1}}\mint_{4B_{\iota}}H(x,Du) \ \dx\nonumber \\
&\stackrel{\eqref{cz3}}{\le}&c_{\texttt{d}}J^{2\bar{p}\Gamma}\rr_{\iota}^{\tau_{0}}\nr{H(\cdot,Du)}_{L^{1+\delta_{\gamma}}(4B_{\iota})}^{\frac{p_{\nu}-p}{p}}\mint_{4B_{\iota}}H(x,Du) \ \dx+\frac{c_{\texttt{nd}}}{J^{p-1}}\mint_{4B_{\iota}}H(x,Du) \ \dx\nonumber \\
&+&c_{\texttt{d}}J^{2\bar{p}\Gamma}\rr_{\iota}^{\tau_{0}}\left(\mint_{4B_{\iota}}H(x,Du) \ \dx\right)\sum_{m=1}^{\kk}\nr{H(\cdot,Du)}_{L^{1+\delta_{\gamma}}(4B_{\iota})}^{\frac{p_{m}-p}{p}}\nonumber \\
&\le&\left(c_{\texttt{d}}J^{2\bar{p}\Gamma}\rr_{\iota}^{\tau_{0}}+\frac{c_{\texttt{nd}}}{J^{p-1}}\right)\mint_{4B_{\iota}}H(x,Du) \ \dx.
\end{eqnarray*}
Summing the content of the previous display over $\nu\in \texttt{I}_{\kk}$ we obtain
\begin{eqnarray}\label{cz11}
\sum_{\nu=1}^{\kk}\mbox{(I)}_{\nu}&\le&\left(c_{\texttt{d}}J^{2\bar{p}\Gamma}\rr_{\iota}^{\tau_{0}}+\frac{c_{\texttt{nd}}}{J^{p-1}}\right)\mint_{4B_{\iota}}H(x,Du) \ \dx.
\end{eqnarray}

\subsubsection*{Step 10 - Matching phases and comparison estimates} Combining \eqref{cz5}, \eqref{cz7}, \eqref{cz10} and \eqref{cz11} we obtain
\begin{eqnarray}\label{cz13}
\mathcal{V}_{0}(Dv_{\iota},Dw_{\iota})&\le&c\left(c_{\texttt{m}}J\rr_{\iota}^{\sigma_{0}}+c_{\texttt{d}}J^{2\bar{p}\Gamma}\rr_{\iota}^{\tau_{0}}+\frac{c_{\texttt{nd}}}{J^{p-1}}\right)\mint_{4B_{\iota}}H(x,Du) \ \dx,
\end{eqnarray}
with $c\equiv c(n,\Lambda,p,p_{1},\cdots,p_{\kk},\kk)$, so via triangular inequality we get
\begin{eqnarray}\label{cz12}
\mathcal{V}(Du,Dw_{\iota})&\le&c\mathcal{V}_{0}(Dv_{\iota},Dw_{\iota})+c\mathcal{V}(Du,Dv_{\iota})\nonumber \\
&\stackrel{\eqref{cz4},\eqref{cz13}}{\le}&c\left(\varepsilon+c_{\texttt{m}}Jr^{\sigma_{0}}+c_{\texttt{d}}J^{2\bar{p}\Gamma}r^{\tau_{0}}+\frac{c_{\texttt{nd}}}{J^{p-1}}\right)\mint_{4B_{\iota}}H(x,Du) \ \dx\nonumber \\
&+&\frac{c}{\varepsilon^{\bar{p}}}\mint_{4B_{\iota}}H(x,F) \ \dx,
\end{eqnarray}
for $c\equiv c(n,\Lambda,p,p_{1},\cdots,p_{\kk},\kk)$. Here we also used that $\rr_{\iota}\le r$, cf. \emph{Step 2}. Next, we set
\begin{eqnarray*}
\mathcal{S}(\varepsilon,r,J,M):=c\varepsilon+cc_{\texttt{m}}Jr^{\sigma_{0}}+cc_{\texttt{d}}J^{2\bar{p}\Gamma}r^{\tau_{0}}+\frac{cc_{\texttt{nd}}}{J^{p-1}}+\frac{c}{\varepsilon^{\bar{p}}M},
\end{eqnarray*}
with $c\equiv c(n,\Lambda,p,p_{1},\cdots,p_{\kk},\kk)$ and use the informations contained in \cite[(4.14)$_{2}$]{comicz} (which come from a covering and exit time argument, so they do not depend on the particular form of $H(\cdot)$ therefore apply in our case as well) to establish that
\begin{eqnarray}\label{cz14}
\mathcal{V}(Du,Dw_{\iota})\le \mathcal{S}(\varepsilon,r,J,M)\lambda,
\end{eqnarray}
which holds for any $J\ge 4$ and for all balls $B_{\iota}$ from the covering in \emph{Step 2}. We stress that \eqref{cz14} holds true independently from the degenerate/nondegenerate/mixed status of $H(\cdot)$. Next, we show that
\begin{eqnarray}\label{cz15}
\mint_{2B_{\iota}}H_{+}(Dw_{\iota}) \ \dx\le c\lambda,
\end{eqnarray}
with $c\equiv c(\texttt{data}_{\textnormal{cz}})$. Assume first that $\texttt{ndeg}_{\textnormal{J}}(2B_{\iota})$ holds with $J=10$. Then we have
\begin{eqnarray*}
\mint_{2B_{\iota}}H_{+}(Dw_{\iota}) \ \dx&\stackrel{\eqref{el2}_{2}}{\le}&\mint_{2B_{\iota}}H_{+}(Dv_{\iota}) \ \dx\nonumber \\
&\stackrel{\eqref{cz6}}{\le}&c\mint_{2B_{\iota}}H(x,Dv_{\iota}) \ dx\le c\mint_{4B_{\iota}}H(x,Du) \ \dx\le c\lambda,
\end{eqnarray*}
with $c\equiv c(\texttt{data})$. On the other hand, if $\texttt{deg}_{\textnormal{J}}(2B_{\iota})$ or $\texttt{mix}_{\textnormal{J}}(2B_{\iota})$ hold again with $J=10$, we have
\begin{eqnarray*}
\mint_{2B_{\iota}}H_{+}(Dw_{\iota}) \ \dx&\stackrel{\eqref{el2}_{2}}{\le}&\mint_{2B_{\iota}}H_{+}(Dv_{\iota}) \ \dx\nonumber \\
&\stackrel{\eqref{cz3}_{1}}{\le}&c\sum_{\nu=1}^{\kk}\rr_{\iota}^{\alpha_{\nu}}\mint_{2B_{\iota}}\snr{Dv_{\iota}}^{p_{\nu}} \ \dx+c\mint_{4B_{\iota}}H(x,Du) \ \dx\nonumber \\
&\stackrel{\eqref{hf}}{\le}&c\sum_{\nu=1}^{\kk}\rr_{\iota}^{\alpha_{\nu}-p_{\nu}\mu}\left(\mint_{4B_{\iota}}H(x,Dv_{\iota})^{1+\sigma_{g}} \ \dx\right)^{\frac{p_{\nu}-p}{p(1+\sigma_{g})}}\mint_{4B_{\iota}}H(x,Dv_{\iota}) \ \dx\nonumber \\
&+&c\mint_{4B_{\iota}}H(x,Du) \ \dx\nonumber\\
&\stackrel{\eqref{cz3}}{\le}&c\rr^{\sigma_{0}}_{\iota}\left(\mint_{4B_{\iota}}H(x,Du) \ \dx\right)\sum_{\nu=1}^{\kk}\nr{H(\cdot,Du)}_{L^{1+\delta_{\gamma}}(4B_{\iota})}^{\frac{p_{\nu}-p}{p}}+c\mint_{4B_{\iota}}H(x,Du) \ \dx\nonumber \\
&\le&c\mint_{4B_{\iota}}H(x,Du) \ \dx \le c\lambda,
\end{eqnarray*}
for $c\equiv c(\texttt{data}_{\textnormal{cz}})$ and \eqref{cz15} is completely proven.
\subsubsection*{Step 11 - A priori estimates for $Dw_{\iota}$} Notice that the frozen integrands $H_{+}(\cdot)$ falls into the realm of those treated in \cite{li2}; in particular it is
\begin{eqnarray*}
\sup_{x\in B_{\iota}}H_{+}(Dw_{\iota})\le c\mint_{2B_{\iota}}H_{+}(Dw_{\iota}) \ \dx \stackrel{\eqref{cz15}}{\le}c\lambda \ \Longrightarrow \ \sup_{x\in B_{\iota}}H(\cdot,Dw_{\iota})\le c_{*}\lambda,
\end{eqnarray*}
with $c,c_{*}\equiv c,c_{*}(\texttt{data}_{\textnormal{cz}})$, where we used the definition in \eqref{++}. At this stage, we can proceed exactly as in \cite[Steps 10-11]{comicz} to first determine $J\equiv J(\texttt{data}_{\textnormal{cz}})\ge 4$, then $\varepsilon\equiv \varepsilon(\texttt{data}_{\textnormal{cz}})\in (0,1)$, $M\equiv M(\texttt{data}_{\textnormal{cz}})$ and finally the threshold radius $r_{*}\equiv r_{*}(\texttt{data}_{\textnormal{cz}})\in (0,1]$ to obtain \eqref{czz} and the proof is complete.

\end{document}